\documentclass[aop,MSNbibl,citesort,dvips]{arximspdf}

\doi{10.1214/10-AOP614}
\volume{40}
\issue{1}
\pubyear{2012}
\firstpage{372}
\lastpage{400}

\makeatletter

\newcommand{\fracc}[2]{{#1/(#2)}}

\renewcommand{\epsilon}{\varepsilon}
\newcommand{\cal}{\mathcal}
\newcommand{\re}{\mathbb{R}}
\newcommand{\p}{\mathbb{P}}
\newcommand{\e}{\mathbb{E}}
\newcommand{\F}{\mathcal F}
\newcommand{\vat}{\mathop{\longrightarrow}_{t\to\infty}}
\newcommand{\ov}{\overline}
\def\i{\infty}
\def\ovl{\overline}
\def\1{\mathbf{1}}

\newcommand{\Uposset}{\mathcal U_m}
\newcommand{\Uinvset}{\mathcal U_{m}^+}
\newproclaim{remark}{Remark}
\newcommand{\eqref}[1]{(\ref{#1})}
\newtheorem{theorem}{Theorem}[section]
\newtheorem{prop}[theorem]{Proposition}
\newproclaim{definition}[theorem]{Definition}
\newtheorem{lemma}[theorem]{Lemma}
\newtheorem{coro}[theorem]{Corollary}
\newproclaim{rem}[theorem]{Remark}
\makeatother

\begin{document}
\begin{frontmatter}

\title{On Az\'{e}ma--Yor processes, their optimal properties and the
Bachelier--drawdown equation}
\runtitle{Azema--Yor processes and Bachelier--Drawdown equation}
\pdftitle{On Azema--Yor processes, their optimal properties and the
Bachelier--drawdown equation}

\begin{aug}
\author[A]{\fnms{Laurent} \snm{Carraro}\ead[label=e1]{laurent.carraro@telecom-st-etienne.fr}},
\author[B]{\fnms{Nicole} \snm{El Karoui}\thanksref{t2}\ead[label=e2]{nicole.elkaroui@cmap.polytechnique.fr}}
\and
\author[C]{\fnms{Jan} \snm{Ob\l\'{o}j}\corref{}\thanksref{t3}\ead[label=e3]{obloj@maths.ox.ac.uk}\ead[label=u1,url]{http://www.maths.ox.ac.uk/\textasciitilde obloj}}
\runauthor{L. Carraro, N. El Karoui and J. Ob\l\'{o}j}
\affiliation{Universit\'{e} de Lyon, Universit\'{e} Paris VI and
University of Oxford}
\address[A]{L. Carraro\\
Telecom Saint-Etienne\\
Universit\'{e} de Lyon\\
25 Rue du Docteur R\'{e}my Annino\\
42000 Saint-Etienne\\
France\\
\printead{e1}} 
\address[B]{N. El Karoui\\
LPMA, UMR 7599\\
Universit\'{e} Paris VI\\
BC 188, 4 Place Jussieu\\
75252 Paris Cedex 05\\
France\\ \printead{e2}}
\address[C]{J. Ob\l\'{o}j\\
Mathematical Institute\\
University of Oxford\\
Oxford OX1 3LB\\
United Kingdom\\ \printead{e3}\\\printead{u1}}
\end{aug}
\pdfauthor{Laurent Carraro, Nicole El Karoui, Jan Obloj}
\thankstext{t2}{Supported by the ``Chaire Risques Financiers'' of the Risk
Foundation, Paris.}
\thankstext{t3}{Supported in part by a Marie Curie
Intra-European Fellowship at Imperial College London within the 6{th}
European Community Framework Programme and the Oxford-Man Institute of Quantitative Finance.}

\received{\smonth{9} \syear{2009}}
\revised{\smonth{5} \syear{2010}}

%
\begin{abstract}
We study the class of Az\'{e}ma--Yor processes defined from a general
semimartingale with a
continuous running maximum process. We show that they arise as unique
strong solutions of the
Bachelier stochastic differential equation which we prove is equivalent
to the drawdown equation.
Solutions of the latter have the drawdown property: they always stay
above a given function of their
past maximum. We then show that any process which satisfies the
drawdown property is in fact an Az\'{e}ma--Yor process.
The proofs exploit group structure
of the set of Az\'{e}ma--Yor processes, indexed by functions, which we
introduce.

We investigate in detail Az\'{e}ma--Yor martingales defined from
a~nonnegative local martingale converging to zero at infinity.
We establish relations between average value at risk, drawdown
function, Hardy--Littlewood transform and its inverse. In particular,
we construct Az\'{e}ma--Yor martingales with a given terminal law and
this allows us to rediscover the Az\'{e}ma--Yor solution to the
Skorokhod embedding problem. Finally, we characterize Az\'{e}ma--Yor
martingales showing they are optimal relative to the concave
ordering of terminal variables among martingales whose maximum
dominates stochastically a given benchmark.
\end{abstract}

%
\pdfkeywords{60G44, 60H10, Azema--Yor process, Bachelier--drawdown equation, drawdown, Hardy--Littlewood transform, average value at risk,
concave order, stochastic order, Skorokhod embedding problem}
\end{frontmatter}

\section{Introduction}\label{intro}

In~\cite{ay1} Az\'{e}ma and Yor introduced a family of simple local
martingales,
associated with Brownian motion or more generally with a~continuous
martingale, which they exploited to solve the Skorokhod embedding
problem. These\vadjust{\goodbreak} processes, called Az\'{e}ma--Yor processes, are simply
functions of the underlying process $X$ and its running maximum $\ovl
X_t=\sup_{s\leq t}X_s$. They proved to be very useful especially in
describing laws of the maximum or of the last passage times of a
martingale and were applied in problems ranging from Skorokhod
embeddings, through optimal inequalities, to Brownian penalizations
(cf. Az\'{e}ma and Yor~\cite{ay1,ay2}, Ob\l\'{o}j and Yor~\cite{oy},
Roynette, Vallois and Yor~\cite{RVY06}). The appearance of Az\'{e}ma--Yor martingales
in all these problems was partially explained
with a characterization in Ob\l\'{o}j~\cite{jamartcarac} as the only
local martingales which can be written as a~function of the couple
$(X,\ovl X_t)$.

Recently these processes have seen a revived interest with applications
in mathematical finance including re-interpretation of classical pricing
formulae (see Madan, Roynette and Yor~\cite{MRY08}) and portfolio
optimization under pathwise constraints (see El Karoui and Meziou \cite
{ELKM06,ELKM08}).
In this paper we uncover a more general structure of these processes
and present new characterizations. We explore in depth their properties
and present some further applications of Az\'{e}ma--Yor processes.
We work in a general setup and extend the concept of Az\'{e}ma--Yor
processes $M^U(X)$, as defined in \eqref{eq:aydevelop} below, to the
context of an arbitrary semimartingale $(X_t)$ with a continuous
running maximum process $\ovl X_t$.

We start by studying the set of Az\'{e}ma--Yor processes $M^U(X)$,
indexed by increasing absolutely continuous functions $U$, and show
that it has a~simple group structure. This allows us to see any
semimartingale with continuous running maximum as an Az\'{e}ma--Yor
process. The main contribution of Section~\ref{sec:bach} is to study
how such representations arise naturally. We show that Az\'{e}ma--Yor
processes allow us to solve explicitly the Bachelier equation, which we
also identify with the drawdown equation. The solutions to the latter
satisfy the drawdown constraint $Y_t\geq w(\ovl Y_t)$. Conversely, if
$(Y_t)$ satisfies the drawdown constraint up to time $\zeta$, then it
can be written as $M^U_{t\land\zeta}(X)$
for some nonnegative $X$. Further, if $Y_\zeta=w(\ovl Y_\zeta)$
a.s., then the inverse process $(X_{t\land\zeta})$ is stopped upon
existing an interval $(0,b)$. We provide explicit relation between
function $U$ which generates Az\'{e}ma--Yor process and functions $w$
and $\varphi$ which feature in the drawdown constraint and in the
SDEs. This characterizes the processes both in a pathwise manner and
differential manner.

Then in Section~\ref{sec:max_HL_AVaR} we specialize further and
investigate Az\'{e}ma--Yor processes
defined from $X=N$ a nonnegative local martingale with continuous
maximum process and with $N_t\to0$ as $t\to\i$. We show how one can
identify explicitly Az\'{e}ma--Yor processes from their terminal
values. In Section~\ref{sec:AVaR_HL} we discuss the average value at
risk and the Hardy--Littlewood transform in a~unified manner using tail
quantiles of probability measures. Then we construct Az\'{e}ma--Yor
martingales with a prescribed terminal law. This allows us to
rediscover, in Section~\ref{sec:skoro}, the Az\'{e}ma--Yor~\cite{ay1}
solution to the Skorokhod embedding problem and give it a new interpretation.

Finally, in the last section, we apply the previous results to uncover
optimal properties of Az\'{e}ma--Yor martingales. More precisely, we
show that all uniformly integrable martingales whose maximum dominates
stochastically a given floor distribution are dominated by an Az\'
{e}ma--Yor martingale in the concave ordering of terminal values.
This problem is an extension of the more intuitive problem, motivated
by finance, to find an optimal martingale for the concave order
dominating (pathwise) a given floor process. It is rather surprising to
find that the two problems have the same solution.
We recover in this way the $\Delta$ operator of Kertz and R\"{o}sler
\cite{KertzRosler92b} and give a~direct way to compute it. These dual
results are compared with the classical primal result stating that
among all uniformly integrable martingales with a fixed terminal law
the Az\'{e}ma--Yor martingale has the largest maximum (relative to the
stochastic order). Furthermore, in both problems we can show that any
optimal martingale is necessarily an appropriate Az\'{e}ma--Yor martingale.


\section{The set of Az\'{e}ma--Yor processes}
\label{sec:def}
Throughout, all processes are defined on $(\Omega,
\mathcal{F}, (\mathcal F_t),P)$ a
filtered probability space satisfying the
usual hypothesis and assumed to be taken right-continuous
with left limits (c\`{a}dl\`{a}g), up to $\i$ included if needed.
All functions are assumed to be Borel measurable.
Given a process $(X_t)$ we denote its running maximum $\ovl
X_t=\sup_{s\le t}X_s$. In what follows, we are essentially
concerned with semimartingales with continuous running maximum,
that we call \textit{max-continuous semimartingales}. Observe that
under this assumption, the process $\ovl X_t=\sup_{s\le t}X_s$
only increases when $\ovl X_t=X_t$ or equivalently
%
%
\begin{equation}
\label{eq:max}
\int_0^T (\ovl X_t-X_t)\,d\ovl X_t=0 \qquad\mbox{for any }T>0.
\end{equation}
%
We let $\tau^b(X)=\tau^b_X=\inf\{t\geq0\dvtx  X_t \geq b\}$
be the first up-crossing time of the level $b$ by process $X$, with the
standard convention that $\inf\{\varnothing\}=\i$.
Note that by max-continuity $X_{\tau^b(X)}=b,$ if $0<\tau^b(X)<\infty
$. With a slight abuse of notation, $\tau^\i_X$
denotes the explosion time of $X$.

\subsection{Definition and properties}
There are two different ways to introduce Az\'{e}ma--Yor processes,
and their equivalence has been proven by several authors (see the
comments below).
%
\begin{definition}\label{def:ay}
Let $(X_t)$ be a max-continuous semimartingale starting from $X_0=a$,
and $\ovl X_t$ its (continuous) running maximum.

With any locally bounded Borel function $u$ we associate the primitive
$U$, defined on $[a,+\infty)$, with initial condition $a^*$, that
is,  $U(x)=a^*+\int_a^x u(s)\,ds$.
The Az\'{e}ma--Yor\vadjust{\goodbreak} process associated with $U$ and $X$ is defined by
one of these two
equations:
%
%
\begin{eqnarray}
\label{eq:aydevelop}
M^U_t(X)
&:=&
U(\ovl X_t)-u(\ovl X_t)(\ovl
X_t-X_t) \quad \mbox{or}\\\label
{eq:aystochinteg}
&\;=&
a^*+\int_0^t u(\ovl X_s)\,dX_s.
\end{eqnarray}
In consequence, $M^U(X)$ is a semimartingale and it is a local
martingale when $X$ is a local martingale.
\end{definition}

Observe that the process $M^U(X)$ is c\`{a}dl\`{a}g, since $U(\ovl X)$
is continuous and
$u(\ovl X_t)(\ovl X_t-X_t)$ is nonzero only on the intervals of
constancy of $\ovl X_t$, where the nonregular process $u(\ovl X_t)$ is
constant. Moreover the jumps of $M^U_t(X)$ are given explicitly by
$-u(\ovl X_t)(X_{t}-X_{t-})$.

We note also that, when $u$ is defined only on some interval $[a,b)$
but $U(b)$ is well defined (finite or infinite), then
we can still define $M^U_t(X)$ for $t\leq\tau^b(X)$ and $M^U_{\tau
^b(X)}(X)=U(b)$. Further, using regularity of paths of $(X_t)$ we have
that \eqref{eq:aydevelop}--\eqref{eq:aystochinteg} hold with $t\land
\tau^b(X)$ instead of $t$ and $u(b):=0$.

The symbol $M^U(X)$ is a slight abuse of notation since this process
depends explicitly on the derivative $u$ rather than the function $U$.
The equivalence between \eqref{eq:aydevelop} and \eqref
{eq:aystochinteg} is easy to establish when $u$ is smooth enough to
apply It\^{o}'s formula,
since the continuity of the running maximum implies from \eqref
{eq:max} that $\int_0^t(\ovl X_t-X_t)\,du(\ovl X_t)\equiv0$.
These results may be extended to all bounded functions $u$ via monotone
class theorem and to all locally bounded functions $u$ via a
localization argument.
Alternatively, the equivalence can be argued using the general balayage
formula (see Nikeghbali and
Yor~\cite{NY06}). The case of locally integrable function $u$ can be
attained for a continuous local martingale $X$, as shown in Ob\l\'{o}j
and Yor~\cite{oy}.

The family of processes in \eqref{eq:aydevelop}--\eqref
{eq:aystochinteg}, or their analog with the local time in zero
replacing the running maximum, were first exhibited when $(X_t)$ is a
continuous local martingale. Az\'{e}ma~\cite{Azema78} obtained them
as an application of formulae for dual predictable projections
resulting from supermartingale representations, and Az\'{e}ma and Yor
\cite{AzemaYor78a} as direct application of their balayage formula.
Az\'{e}ma and Yor~\cite{ay1,ay2} described these processes in more
detail and used them to solve the Skorokhod embedding problem.


The importance of the family of Az\'{e}ma--Yor martingales is well
exhibited by Ob\l\'{o}j~\cite{jamartcarac} who proves that in the
case of a continuous local martingale $(X_t)$ all local
martingales which are functions of the couple $(X_t,\ovl X_t)$,
$M_t=H(X_t,\ovl X_t)$ can be represented as a $M=M^U(X)$ local
martingale associated with a locally integrable function $u$. We
note that such processes are sometimes called
\textit{max-martingales}.
\subsection{Monotonic transformations and Az\'{e}ma--Yor processes}
\label{subsection:monotonic}
We want to investigate further the structure of the set
of Az\'{e}ma--Yor processes associated with a max-continuous\vadjust{\goodbreak}
semimartingale $(X_t)$.
One of the most remarkable properties of these processes is that
their running maximum can be easily computed, when the function
$U$ is nondecreasing ($u\geq0$).

We denote by $\Uposset$ the set of such functions, that is
absolutely continuous functions defined on an appropriate interval
with a locally bounded and nonnegative derivative. This set is
stable by composition, that is, if $U$ and~$F$ are in $\Uposset$,
and defined on appropriate intervals, then $U\circ F(x)=U(F(x))$ is
in $\Uposset$. We let $\Uinvset$ be the set of increasing
functions $U\in\Uposset$, with inverse function $V \in\Uposset$,
or equivalently of functions $U$ such that $u=U'>0$, and both $u$
and $1/u$ are locally bounded. Throughout, when we consider an
inverse function $V$ of $U\in\Uinvset$ then we choose $v(y)=V'(y)=1/u(V(y))$.

In light of \eqref{eq:aydevelop}, we then have:
%
\begin{prop}\label{prop:group}
\textup{(a)} Let $U\in\Uposset$, $X$ be a max-continuous semimartingale and $M^U(X)$
be the Az\'{e}ma--Yor process in \eqref{eq:aydevelop}.
Then
%
%
\begin{equation}
\label{eq:maxAY}
\ovl{M^U_t(X)}=U(\ovl X_t),
\end{equation}
and $M^U(X)$ is a max-continuous semimartingale.

\textup{(b)} Let $F\in\Uposset$ defined on an appropriate interval, so that
$U\circ F$ is well defined. Then
\[
M^U_t (M^F(X) )=M^{U\circ F}_t(X).
\]
\end{prop}

%
\begin{rem}
It follows from point (b) above that the set of Az\'{e}ma--Yor processes
indexed by $U \in
\Uinvset$ defined on whole $\re$ with $U(\re)=\re$, is a group
under the operation $\otimes$ defined by
\[
M^U\otimes
M^F:=M^{U\circ F}.
\]
Note that $M^{\operatorname{Id}}_t(X)=X_t$, where $\operatorname{Id}(x)=x$ is the
identity mapping.
\end{rem}

\begin{pf*}{Proof of Proposition~\ref{prop:group}}
(a) In light of \eqref{eq:aydevelop}, when $u$ is nonnegative, the
Az\'{e}ma--Yor process $M^U_t(X)$ is dominated by $U(\ovl X_t)$,
with equality if $t$ is a point of increase of $ \ovl X_t$. Since
$U$ is nondecreasing we obtain \eqref{eq:maxAY}. Moreover, since
$U(\ovl X)$ is a continuous process, $M^U(X)$ is a max-continuous
semimartingale and we may take an Az\'{e}ma--Yor process of it.

(b) Let $F$ be in $\Uposset$, $f=F'\geq0$, such that
$U\circ F$ is well defined. We have from \eqref{eq:maxAY} and \eqref
{eq:aydevelop}
%
%
\begin{eqnarray}
M^U_t (M^F(X))
&=&U (F(\ov X_t) )-u (F(\ov X_t)
)f(\ov X_t) (\ov X_t-X_t)\nonumber
\\[-8pt]
\\[-8pt]
&=&M^{U\circ F}_t(X),
\nonumber
\end{eqnarray}
where we used $ (U(F(x)) )'=u(F(x))f(x)$.
\end{pf*}

The two properties described in Proposition~\ref{prop:group} are rather
simple but extremely useful.
We phrase part (b) above for stopped processes and for $F=V=U^{-1}$ as a
separate corollary.

%
\begin{coro}\label{cor:dual}
Let $a<b\le\infty$, $U\in\Uinvset$ the primitive function of
a~locally bounded $u\dvtx [a,b)\to(0,\i)$
with $U(a)=a^*$. Let $V\dvtx [a^*,U(b))\to[a,b)$ be the inverse of $U$
with locally bounded derivative $v(y)=1/u(V(y))$.

Then for any max-continuous semimartingale $(X_t)$, $X_0=a$, stopped at
the time
$\tau^b=\tau^b(X)=\inf\{t\dvtx X_t\geq b\}$ we have
%
%
\begin{equation}
\label{eq:inverse}
X_{t\land\tau^b}=M^V_{t\land\tau^b}(M^U(X)).
\end{equation}
From the differential point of view, on $[0,\tau^b)$,
%
%
\begin{equation}
\label{eq:inversedif}
dY_t=u (\ovl X_t )\,dX_t\quad\mbox{and}\quad dX_t=v
(\ovl Y_t )\,dY_t \qquad\mbox{where }Y_t=M^U_t(X).
\end{equation}
\end{coro}

Consider $u$ as above with $b=U(b)=\infty$. As a consequence of
the above, any max-continuous semimartingale $(X_t)$ can be seen
as an Az\'{e}ma--Yor process associated with $U$. Indeed,
$X_t=M^U_t(Y)$ with $Y_t=M^V_t(X)$. In the following section we
study how such representations arise in a natural way.

\section{The Bachelier--drawdown equation}
\label{sec:bach}
In his paper, ``Th\'{e}orie des probabilit\'{e}s continues,'' published
in 1906,
French mathematician Louis Bachelier~\cite{Bach} was the first to consider
and study stochastic differential equations. Obviously,
he did not prove in his paper existence and uniqueness results but focused
his attention on some particular types of SDEs. In this way,
he obtained the general structure of processes with independent
increments and continuous paths, the definition of diffusions
(in particular, he solved the Langevin equation),
and generalized these concepts to higher dimensions.
%
\subsection{The Bachelier equation} \label{subsection:Bachequation}
In particular, Bachelier (\cite{Bach}, pages 287--290) considered and
``solved'' an
SDE depending on the maximum of the
solution, $dY_t=\varphi(\ovl Y_t)\,dX_t$ which we call the
Bachelier equation.
Let $U\in\Uinvset$ and $V\in\Uinvset$ its inverse function with
derivative $v$.
From \eqref{eq:inversedif} and \eqref{eq:maxAY} we see that the Az\'{e}ma--Yor process $Y=M^U(X)$ verifies
the Bachelier equation for $\varphi(y)=1/v(y)$. Now, we can solve the
Bachelier equation
as an inverse problem. We present a rigorous, simple and explicit
solution to this equation. We note that a similar approach is developed
in Revuz and Yor~\cite{ry}, Exercice~VI.4.21.
%
\begin{theorem}\label{thm:bach}
Let $(X_t\dvtx t\ge0)$, $X_0=a$, be a max-continuous semimartingale.
Consider a positive Borel function $\varphi\dvtx [a^*,\i)\to(0,\i)$ such
that $\varphi$ and $1/\varphi$ are locally bounded. Let
$V(y)=a+\int_{a^*}^y\frac{ds}{\varphi(s)}$, and $U$ its inverse defined on $(a,V(\i))$.

The Bachelier equation,
%
%
\begin{equation}\label{eq:bach}
dY_t=\varphi(\ovl Y_t)\,dX_t, \qquad Y_0=a^*,
\end{equation}
has a strong, pathwise unique, max-continuous solution defined up
to its explosion time $\tau^\i_Y=\tau^{V(\i)}_X$
given by $Y_t=M^U_t(X),  t< \tau^{V(\i)}_X$.\vadjust{\goodbreak}

 When $X$ is a
continuous local martingale it suffices to assume that $1/\varphi$ is a
locally integrable function.
\end{theorem}

%
%
\begin{pf}
The assumptions on $\varphi$ imply that $V$ and therefore $U$ are
in~$\Uinvset$ with $U(a)=a^*$. With $u=\varphi(V)$, Definition \ref
{def:ay} and \eqref{eq:maxAY} gives that the Az\'{e}ma--Yor process
$M^U(X)$ verifies
\[
dM^U_t(X)= u(\ovl X_t)\,dX_t = \varphi(\ovl{M^U_t(X)})\,dX_t,\qquad t<\tau^{V(\i)}_X.
\]
Furthermore, on $\tau^{V(\i)}(X)<\infty,$ $M^U_{\tau
^{V(n)}}(X)=U(V(n))=n$ and we see that if $V(\i)<\i$ then
$\tau^{V(\i)}_X$ is the explosion time of $M^U(X)$. So,
$M^U(X)$ is a~solution of \eqref{eq:bach}.

Now let $Y$ be a max-continuous solution to equation
\eqref{eq:bach}. Equation \eqref{eq:aystochinteg} in Definition \ref
{def:ay} and \eqref{eq:bach} imply
that $dM^V_t(Y)=dX_t $ on $[0,\tau^\i_Y)$. It follows from
Corollary~\ref{cor:dual} that $Y_t=M^U_t(X)$ and $\tau^{V(\i)}_X$ is
the
explosion time $\tau^\i_Y$ of $Y$.

The above result extends to more general $\varphi$ whenever $U,V$ and
$M^U(X)$ are well defined. When $X$ is a continuous local martingale,
to define $V$ and~$U$ it is sufficient (and necessary) to assume
$1/\varphi$ is locally integrable.
That $M^U(X)$ is then well defined follows from Ob\l\'{o}j and Yor
\cite{oy}.
\end{pf}

The above extends naturally to the
case when $a$ and $a^*$ are some $\mathcal{F}_0$-measurable random
variables. It suffices to assume that $\varphi$ is well defined on
$[l,\infty)$ where $-\i\le l$ is the lower bound of the support
of $a^*$.
We could also consider $X$ which is only defined up to its
explosion time $\tau^\i_X$ which would induce $\tau^\i_Y=\tau^\i
_X\land
\tau_X^{V(\i)}$.

In Section~\ref{sec:max_HL_AVaR} we will also consider the case when
$\varphi\equiv0$ on $(r,\infty)$ and then $(Y_t)$ is stopped upon
hitting $r$.

Finally note that under a stronger assumption that $X$
has no positive jumps, \textit{any} solution of the Bachelier
equation has no positive jumps and hence is a max-continuous semimartingale.

\subsection{Drawdown constraint and drawdown equation}
\label{sec:dd}
In various applications, in particular in financial mathematics,
one is interested in processes which remain above a (given)
function $w$ of their running maximum. The purpose of this section
is to show that Az\'{e}ma--Yor processes provide a direct answer to
this problem when the underlying process $X$ is positive.
The following notion will be central throughout the rest of the paper.
%
\begin{definition}\label{def:DD}
Given a function $w$, we say that a c\`{a}dl\`{a}g process $(M_t)$
satisfies $w$-\textit{drawdown}
($w$-DD) constraint up to the (stopping) time $\zeta$, if $\min\{
M_{t-}, M_t\}> w(\ovl M_t)$ for all $0\leq t<\zeta$ a.s.
\end{definition}

 We will see in Section~\ref{sec:max_HL_AVaR} that for a
local martingale $M$ it suffices
to impose $M_t>w(\ovl M_t)$ in the above definition.\vadjust{\goodbreak}

Az\'{e}ma--Yor processes, $Y=M^U(X)$ defined from a \textit{positive}
max-continuous semimartingale $X$ and function $U\in\Uinvset$
provide an example of such processes with DD-constraint function
$w$ defined from $U$ and $V=U^{-1}$ by
%
%
\begin{equation}\label{eq:DDfunction}
w(y)=h(V(y))=y-V(y)/v(y) \qquad\mbox{where }  h(x)=U(x)-x u(x).
\end{equation}
Indeed, thanks to the positivity of $X$ and $u$ we have
%
%
\begin{eqnarray}\label{eq:DDforAY}
\min\{Y_{t-},Y_t\}&=&U(\ovl X_t)-u(\ovl X_t)\ovl X_t+u(\ovl X_t)\cdot
\min\{X_{t-},X_t\}\nonumber
\\[-9pt]
\\[-9pt]
&>& U(\ovl X_t)-u(\ovl X_t)\ovl X_t=h(\ovl X_t)=w(\ovl Y_t).
\nonumber
\end{eqnarray}
The converse is more interesting. We show below that if we start with
a~given $w$ then $M^U(X)$ satisfies the $w$-DD constraint for $U=V^{-1}$
and~$V$ given in~\eqref{eq:DDEquation} below. Furthermore, it turns
out that all processes which satisfy a~drawdown constraint are of this
type. More precisely,
given a max-continuous semimartingale $Y$ satisfying the $w$-DD
constraint we can find explicitly $X$ such that $Y$ is the
Az\'{e}ma--Yor process $M^U(X)$.
Moreover, the first
instant~$Y$ violates the drawdown constraint is precisely the
first hitting time to zero of $X$. For a precise statement we need
to introduce the set of admissible functions $w$.\vspace*{-2pt}
%
\begin{definition}\label{def:w-DD}
We say that $w\dvtx [a^*,\i]\to\re$ is a \textit{drawdown function} if it
is nondecreasing and there exists $r_w\leq\i$ such that $y-w(y)>0$
is locally bounded and locally
bounded away from zero on $[a^*,r_w)$ and $w(y)=y$ for $y\geq
r_w$.\vspace*{-2pt}
\end{definition}

We impose $w$ nondecreasing as it is intuitive for applications.
It will also arise naturally in Section~\ref{sec:max_HL_AVaR}.
We introduced here $r_w$ as it gives a convenient way to stop the
process upon hitting a given level and again it will be used in Section
\ref{sec:max_HL_AVaR}. If a drawdown function $w$ is defined on
$[a^*,\i)$, then we put $w(\i)=\lim_{y\uparrow\i} w(y)$, and the
above definition requires that $w(\i)=\i$.
In fact for the results in this section it is not necessary to require
any monotonicity from $w$ or that $w(\infty)=\infty$, we comment this below.

We let $\tau_0(X)=\tau_0^X=\inf\{t\dvtx  \min\{X_{t-},X_t\}\leq0\}$ and
note that when $X$ is nonnegative
then $X_{\tau_0^X}\geq0$ on the set $\{\tau_0^X<\infty\}$. Further
let \label{page:zeta_DD} $\zeta_w^Y=\inf\{t\dvtx\break  \min\{Y_{t-},Y_t\}\leq
w(\ovl Y_t)\}$.
As mentioned before, definitions of both $\tau_0$ and $\zeta_w$
simplify for local martingales (see Lemma~\ref{lem:Nstop}
and Corollary~\ref{cor:N-DD} in Section~\ref{sec:max_HL_AVaR}).\vspace*{-2pt}

\begin{theorem}\label{thm:DD} Consider a drawdown function $w$ of
Definition~\ref{def:w-DD}. Then the solution $V$ of the ODE \eqref
{eq:DDfunction} with $V(a^*)=a>0$ is given as
%
%
\begin{equation}
\label{eq:DDEquation}
V(y)=a\exp \biggl(\int_{a^*}^y\frac
{1}{s-w(s)}\,ds \biggr),\qquad y\geq a^*.
\end{equation}
For a nonnegative max-continuous semimartingale $(X_t)$, $X_0=a$, and
$\zeta:=\tau_0(X)\land\tau^{V(r_w)}(X)$ the drawdown equation
%
%
\begin{equation}\label{eq:bachDD}
dY_t= \bigl(Y_{t_-}-w(\ovl Y_t) \bigr)\frac{dX_t}{X_{t_-}},\qquad
t<\zeta,\vadjust{\goodbreak}
\end{equation}
has a strong, pathwise unique, max-continuous solution which satisfies \mbox{$w$-DD}
constraint on $[0,\zeta)$ and $Y_0=a^*$, given by $Y_t=M^U_t(X)$,
where $U$ is the inverse of $V$.
We have $\zeta_w^Y=\zeta$ and further $Y_{\zeta_w^Y}=w(\ovl Y_{\zeta
_w^Y})$ on $\{X_\zeta\in\{0,V(r_w)\}\}$.

Conversely, given $(Y_t)$ a max-continuous semimartingale satisfying
$w$-DD constraint up to $\zeta:=\zeta_w^Y$, with $Y_0=a^*$, there
exists a pathwise unique max-continuous semimartingale $(X_t\dvtx t< \zeta
)$, $X_0=a$, which solves \eqref{eq:bachDD}. $X$ may be deduced from
$Y$ by the Az\'{e}ma--Yor bijection $X_t=M^V_t(Y)$ and $\zeta=\tau
_0(X)\land\tau^{V(r_w)}(X)$.
\end{theorem}

%
\begin{rem}
Naturally $V(y)\equiv\i$ for $y> r_w$. However, $V(r_w)$ could be
either finite or infinite and consequently $\tau^{V(r_w)}(X)$ can be
either a hitting time of a finite level or the explosion time for $X$.

Observe that $\{X_\zeta\in\{0,V(r_w)\}\}$ could be larger than $\{
\zeta<\i\}$. This will be the case in Section~\ref{sec:max_HL_AVaR}
where $X_t\to0$ as $t\to\i$ and in fact $X_\zeta\in\{0,V(r_w)\}$
a.s. Naturally, we also have $Y_{\zeta_w^Y-}=w(\ovl Y_{\zeta_w^Y})$
on $\{X_{\zeta-}\in\{0,V(r_w)\}\}$. Note also that in the converse
part of the theorem we could have $Y_\zeta<w(\ovl Y_\zeta)$ which
would correspond to $X_\zeta<0$.
\end{rem}

%
\begin{rem}
It will be clear from the proof that the theorem holds without assuming
any monotonicity on $w$ or that $w(\infty)$ is defined and equal to
$\infty$. The only changes are   $Y_{\zeta_w^Y}=w(\ovl Y_{\zeta
_w^Y})$ on $\{X_\zeta=V(r_w)\}$ if and only if $w(r_w)=r_w$ and if
$V(\i)<\i$
then $Y$ explodes at $\tau^{V(\i)}_X$.
\end{rem}

\begin{pf*}{Proof of Theorem~\ref{thm:DD}} Expression for $V$ in
terms of $w$ follows as $v(y)=V(y)/(y-w(y))$. Note that $V(\i)=\i$.
Hence, for $t<\zeta$, $Y_t=M_t^U(X)$ is well defined, and recall from
Corollary~\ref{cor:dual}
that $X_t=M_t^V(Y)$ and $\ovl X_t=V(\ovl Y_t)$. Direct computation yields
\[
Y_{t-}-w(\ovl Y_t)=Y_{t-}-U(\ovl X_t)+u(\ovl X_t)\ovl X_t=u(\ovl
X_t)X_{t-} .
\]
%
Thanks to the positivity of $u$ and $X$ and $X_-$ on $  t<\zeta$, we
have that
$Y_{t-}, Y_t> w(\ovl Y_t)$
and it follows from \eqref{eq:aystochinteg} that $Y=M^U(X)$
solves \eqref{eq:bachDD}.

Now consider any max-continuous solution $Y$ of \eqref{eq:bachDD},
$\min\{Y_{t-},Y_t\}>w(Y_t)$ for $t<\zeta$.
Then, using \eqref{eq:aydevelop} and \eqref{eq:aystochinteg}, we
have
\[
\frac{dY_t}{Y_{t_-}-w(\ovl Y_{t_-})}=
\frac{v(\ovl Y_t)}{M^V_{t_-}(Y)}\,dY_t=
\frac{dM^V_t(Y)}{M^V_{t_-}(Y)} .
\]
Since $Y$ is solution of \eqref{eq:bachDD}, $X$ and $M^V(Y)$ have the
same relative
stochastic differential and the same initial condition.
Then, there are undistinguishable processes and Corollary \ref
{cor:dual} yields $Y_t=M^U_t(X)$.

Finally, when $X_\zeta= 0$ [resp., $X_\zeta=V(r_w)$] we have $Y_\zeta
= U(\ovl X_\zeta)-u(\ovl
X_\zeta)\ovl X_\zeta=w(V(\ovl X_\zeta))=w(\ovl Y_\zeta)$ [resp.,
$Y_\zeta=r_w=\ovl Y_\zeta=w(\ovl Y_\zeta)$] and
$\zeta=\zeta_w^Y$. If $X_\zeta\notin\{0,V(r_w)\}$, then
$X_{\zeta-}=0$ or $\zeta=\i$, and in both cases $\zeta=\zeta_w^Y$.

Consider now the second part of the theorem. We can rewrite \eqref
{eq:bachDD} as
%
%
\begin{equation}\label{eq:DD}
\frac{dY_t}{Y_{t_-}-w(\ovl Y_t)}=\frac{dX
_t}{X_{t_-}},\qquad t< \zeta.
\end{equation}
This equation defines without ambiguity a positive process $X$ starting
from $X_0=a>0$.
By assumption on $w$, the solution $V$ of \eqref{eq:DDfunction} is a positive
finite increasing function on $[a^*,r_w)$, $V(y)/v(y)= y-w(y)$.
Put $\widehat X_t=M^V_t(Y)$, and observe that the differential
properties of $V$ imply that $\widehat X_t=v(\ovl Y_t)(Y_t-w(\ovl
Y_t))>0$, for $t<\zeta$. Then, the stochastic differential of
$M^V_t(Y)=\widehat X_t$
is
\[
d\widehat X_t=v(\ovl Y_t)\,dY_t=\widehat X_{t-}\bigl(Y_{t-}-w(\ovl
Y_t)\bigr)^{-1}\,dY_t,
\]
and hence both $\widehat X$ and $X$ are solutions of the same
stochastic differential equation and have the same initial
conditions. So, they are undistinguishable processes. Identification of
$\zeta$ follows as previously.
\end{pf*}

Naturally, since $Y=M^U(X)$ solves both the Bachelier equation \eqref
{eq:bach} and the drawdown equation \eqref{eq:bachDD} these equations
are equivalent. We phrase this as a corollary in the case $\zeta=\i$
a.s. in \eqref{eq:bachDD}.
%
\begin{coro}\label{coro:DD}
Let $(X_t\dvtx t\ge0)$, $X_0=a$, be a positive max-continuous
semimartingale, $\tau_0^X=\i$ a.s., and $\varphi,V$ as in Theorem
\ref{thm:bach}.
Then, the Bachelier equation \eqref{eq:bach} is equivalent to the
drawdown equation \eqref{eq:bachDD} where $w$ and~$V$ are linked via
\eqref{eq:DDfunction} or equivalently via \eqref{eq:DDEquation}.
\end{coro}

%


The drawdown equation \eqref{eq:bachDD}
was solved previously by Cvitani\'{c} and
Karat\-zas~\cite{CK95} for $w(y)=\gamma y$, $\gamma\in(0,1)$,
and recently by Elie and Touzi~\cite{ET08}.
The use of Az\'{e}ma--Yor processes simplifies considerably the proof
and allows for a general $w$ and $X$.
We have shown that this equation has a unique strong solution and is
equivalent to the Bachelier equation.

Note that we assumed $X$ is positive. The quantity $dX_t/X_{t-}$ has a
natural interpretation as the differential of the stochastic
logarithm of $X$. In various applications, such as financial
mathematics, this logarithm process is often given directly since
$X$ is defined as a stochastic exponential in the first place.

\subsubsection*{An illustrative example}\label{page:ill_ex}Let $X$ be a
positive max-continuous semimartingale
such that $X_0=1$.
Let $U$ be the power utility function defined on~$\re^+$ by
$U(x)=\frac{1}{1-\gamma}  x^{1-\gamma}$ with $0<\gamma< 1$
and $u(x)=x^{-\gamma}$ its derivative. 
The inverse function $V$ of $U$ is $V(y)=((1-\gamma)y)^{1/(1-\gamma)}$
and its derivative is $v(y)=((1-\gamma)y)^{\gamma/(1-\gamma)}$.

Then the (power) Az\'{e}ma--Yor process is
\[
M^U_t(X)=Y_t=\frac{1}{1-\gamma}(\ovl X_t)^{1-\gamma}
\biggl (\gamma+(1-\gamma)\frac{X_t}{\ovl X_t} \biggr)= \ovl
Y_t \biggl(\gamma+(1-\gamma)\frac{X_t}{\ovl X_t} \biggr).
\]
Since $X$ is positive, $Y_t> w(\ovl Y_t)=\gamma\ovl Y_t$. The drawdown
function $w$ is the linear one, $w(y)=\gamma y$.

The process $(Y_t)$ is a semimartingale (local martingale if $X$ is a local
martingale) starting from $Y_0=1$, and
staying in the interval $[\gamma\ovl Y_t,\ovl Y_t]$. Since the power
function $U$ is concave, we also have an other
floor process $Z_t=U(X_t)$. Both processes $Z_t$ and $\gamma\ovl
Y_t=\gamma\ovl Z_t$ are dominated by $Y_t$. They are not comparable in
the sense that in general at time $t$ either one of them can be
greater. We study floor process $Z$ in more detail in Section~\ref{sec:floor}.

The Bachelier--drawdown equation \eqref{eq:bach}--\eqref{eq:bachDD}
becomes here
%
%
\begin{eqnarray}
dY_t&=&{\ovl X}_t^{\,-\gamma}\,dX_t= \bigl((1-\gamma)\ovl Y_t
\bigr)^{-\fracc{\gamma}{1-\gamma}}\,dX_t\nonumber
\\[-8pt]
\\[-8pt]
&=&(Y_{t-}-\gamma\ovl Y_t)\frac{dX_t}{X_{t-}} .
\nonumber
\end{eqnarray}

As noted above, this equation, for
a class of processes $X$, was studied in Cvitani\'{c} and Karatzas
\cite{CK95}. Furthermore, in~\cite{CK95} authors in fact
introduced processes $M^U(X)$ where $U$ is the a power utility
function, and used them to solve the portfolio optimization
problem with drawdown constraint of Grossman and Zhou
\cite{GZ93} (see also~\cite{ET08}). Using our methods we can
simplify and generalize their results and show that the portfolio
optimization problem with drawdown constraint, for a general
utility function and a general drawdown function, is equivalent to
an unconstrained portfolio optimization problem with a modified
utility function. We develop these ideas in a~separate paper.

\section{Setup driven by a nonnegative local martingale converging to
zero}\label{sec:max_HL_AVaR}
In the previous section we investigated Az\'{e}ma--Yor processes built
from a~nonnegative semimartingale as solutions to the drawdown
equation \eqref{eq:bachDD}. We specialize now further and study in
detail Az\'{e}ma--Yor processes associated with $X=N$, a~nonnegative
local martingale converging to zero at infinity.
The maximum of $N$ has a universal law which, together with $N_\i=0$,
allows to write Az\'{e}ma--Yor martingales explicitly from their
terminal values (see Sections~\ref{sec:maximum} and \ref
{sec:ay_given_values}). Our study exploits tail quantiles of
probability measures and is intimately linked with the average value at
risk and the Hardy--Littlewood transform of a measure, as explored in
Section~\ref{sec:AVaR_HL}. Finally, combining these results with
Theorem~\ref{thm:DD}, we construct Az\'{e}ma--Yor martingales with
prescribed terminal distributions and in particular obtain the Az\'
{e}ma--Yor~\cite{ay1} solution of the Skorokhod embedding problem.

\subsection{Universal properties of $X=N$}\label{sec:maximum}
A nonnegative local martingale $(N_t)$ is a supermartingale and it is
a (true) martingale if and only if $\e N_t=\e N_0$, $t\geq0$. We also
have that if $N_t$ or $N_{t-}$ touch zero then $N_t$ remains in zero
(see, e.g., Dellacherie and Meyer~\cite{DellacherieMeyer80}, Theorem~VI.17).
%
\begin{lemma}\label{lem:Nstop}
Consider a nonnegative local martingale $(N_t)$ with $N_{0-}:=N_0>0$. Then
%
%
\begin{equation}\label{eq:Nstop}
\tau_0(N)=\inf\{t\dvtx  N_t=0  \mbox{ or }  N_{t-}=0\}=\inf\{t\dvtx  N_t=0\}
\end{equation}
and $N_u\equiv0$, $u\geq\tau_0(N)$.\vadjust{\goodbreak}
\end{lemma}

This yields an immediate simplification of the $w$-DD condition. In
fact in Definition~\ref{def:DD}, and the definition of $\zeta_w^Y$
before Theorem~\ref{thm:DD}, it suffices to compare $w(\ovl Y_t)$ with
$Y_t$ instead of $Y_t$ \textit{and} $Y_{t-}$.
%
\begin{coro}\label{cor:N-DD}
Let $w$ be a drawdown function of Definition~\ref{def:DD} and
$(Y_t)$ a max-continuous local martingale with $ Y_\zeta=w(\ovl
Y_\zeta)$ a.s. on $\{\zeta<\i\}$, where $\zeta=\inf\{t\dvtx  Y_t\leq
w(\ovl Y_t)\}$. Then $Y$ satisfies $w$-DD condition up to time $\zeta
_w^Y=\zeta$.
\end{coro}

\begin{pf}
Assume $r_w=\i$ and let $N_t=M^V_t(Y)$ where $V$ is given via~\eqref
{eq:DDEquation}. Using \eqref{eq:DDfunction} and \eqref{eq:DDforAY},
similarly as in the proof of Theorem~\ref{thm:DD}, and Definition~\ref
{def:ay}, $(N_t\dvtx t\leq\zeta)$ is a nonnegative max-continuous local
martingale and $\zeta=\inf\{t\dvtx  N_t=0\}$. Using \eqref{eq:Nstop} we
have $\zeta=\tau_0(N)$, and our assumptions also give $N_\zeta=0$ on
$\{\zeta<\i\}$. It follows from Theorem~\ref{thm:DD} that
$Y_t=M^U(X)_t$, $U=V^{-1}$ satisfies the $w$-DD constraint up to $\zeta
$ and $\zeta=\tau_0(N)=\zeta_w^Y$. For the case $r_w<\i$ it
suffices to note that all processes are max-continuous and hence the
first hitting times for $\ovl N_t$ and $\ovl N_{t-}$ are equal.
\end{pf}

Throughout this and following sections, we assume that
%
%
\begin{eqnarray}\label{eq:Ndef}
\begin{tabular}{@{}p{275pt}@{}}
  ($N_t\dvtx  t\geq0)$ \mbox{is a nonnegative max-continuous local
martingale,} $\displaystyle N_t\vat0$\mbox{ a.s.}
\end{tabular}
\\[-10pt]
\nonumber
\end{eqnarray}
%
Note that in particular $\inf\{t\geq0\dvtx  \ovl N_t = \ovl N_\infty\}
<\infty$ a.s.

We recall some well-known results on the distribution of
the maximum of~$N$ (see Exercice III.3.12 in Revuz--Yor
\cite{ry}). We assume that $N_0$ is a~constant. If $N_0$ is random all
results should be read conditionally on $\F_0$.
%
\begin{prop}\label{prop:uniform}
Consider $(N_t)$ satisfying \eqref{eq:Ndef} with $N_0>0$ a constant.
\begin{longlist}[(b)]
\item[(a)] The random variable $N_0/\ovl N_\i$ is uniformly distributed on
$[0,1]$.
\item[(b)] The same result holds for the conditional distribution
in the following sense: let $\ovl N_{t,\i}= \sup_{t\leq u\leq\i}
N_u$ then
\[
\p(K>\ovl N_{t,\i}|\F_t)=(1-N_t/K)^+;
\]
that is,  $\ovl N_{t,\i}$ has the same $\F_t$-conditional
distribution as
$N_t/\xi$ where $\xi$ is an independent uniform variable on $[0,1]$.
\item[(c)] Let $\zeta=\tau_0(N)\wedge\tau^b(N)=\inf\{t\dvtx  N_t\notin(0,b)\}
$, $b>N_0$. $(N_{t\wedge\zeta})$ is a bounded martingale and
$N_{\zeta}\in\{0,b\}$. Furthermore, ${\ovl N}_{\zeta}={\ovl N}_\i
\wedge b$ is
distributed as $(N_0/\xi)\wedge b$, where $\xi$ is uniformly
distributed on $[0,1]$.
\end{longlist}
\end{prop}

\begin{remark*} (a) Given the event $\{\ovl
N_{\zeta} < b\}=\{N_{\zeta}=0\}$, $N_0/\ovl N_{\zeta}$ is uniformly
distributed on $(N_0/b,1]$. The probability of the event $\{\ovl
N_{\zeta} = b\}$ is $N_0/b$.

(b) Any nonnegative martingale $N$ stopped at $\zeta$, with $N_{\zeta
}\in\{0,b\}$ a.s., may be extended into a local martingale (still
denoted by $N$) satisfying~\eqref{eq:Ndef}, by putting $N_t:=N_\zeta
+\1_{\{\ovl N_\zeta=b\}}(N'_t-N'_\zeta),$ $t>\zeta$, where $N'$ is
another local martingale as in~\eqref{eq:Ndef}.
\end{remark*}

\begin{pf*}{Proof of Proposition~\ref{prop:uniform}}
(a) Let us consider the Az\'{e}ma--Yor martingale associated with
$(N_t)$ and the function $U(x) = (K-x)^+$, where $K\geq N_0$ is a fixed
real. Thanks to the positivity of $(N_t)$, the
martingale $M^U(N)$ is bounded by $K$,
\[
0 \leq M^U_t(N) = (K - \ovl N_t)^+ + \1_{\{K > \ovl N_t \}}(\ovl N_t- N_t)=
\1_{\{K > \ovl N_t \}}(K-N_t)\leq K.
\]
So $M^U_t(N)$ is a uniformly integrable martingale, and $\e M^U_\i
(N)=M^U_0(N)$.
Since $N_t\to0$ as $t\to\infty$, $\ovl N_\infty<\infty$ and we
have $M^U_\infty(N)=K\mathbf{1}_{K>\ovl N_\infty}$ and the previous
equality can be written as $K \p(K > {\ovl N}_\i) = K - N_0$, or
equivalently $\p (\frac{N_0}{{\ovl N}_\i} \leq\frac{N_0}{K}
 ) =
\frac{N_0}{K}$. That is exactly the desired result.

(b) This result is the conditional version of the previous one. The
reference process is now the process $(N_{t+h}\dvtx  h\geq0)$
adapted to the filtration $\F_{t+h}$, local martingale for the
conditional probability measure $\p(\cdot|\F_t)$.

(c) From \eqref{eq:Nstop} and since $N$ is nonnegative and
max-continuous it follows that
$\tau_0(N)\wedge\tau^b(N)=\inf\{t\dvtx  N_t\notin(0,b)\}$ and that
$N_{\zeta}=b, \mbox{ or }  0$. Then, we have that $\ovl N_{\zeta
}=\ovl N_\i\wedge b$ since when $N_{\zeta}=b$, the maximum $\ovl
N_{\zeta}$
is also equal to $b$.
\end{pf*}

\subsubsection*{Remark about last passage times}
Recently, for a
continuous local martingale~$N$, Madan, Roynette and Yor
\cite{MRY08} 
have interpreted the event $\{K>\ovl N_{T,\infty}\}$ in terms of the
last passage
time $g_K(N)$ over the level $K$, as $\{K>\ovl
N_{T,\infty}\}=\{g_K(N)<T\}$. Our last proposition yields
immediately their result: the normalized put pay-off is the conditional
probability of $\{g_K(N)<T\} \dvtx(1-N_T/K)^+= \p(g_K(N)<T|\F_T)$. In particular we obtain the whole
dynamics of the put prices,
\[
\e [(K-N_T)^+|\F_t ]= K\p\bigl(g_K(N)<T|\F_t\bigr),\qquad t\leq T,
\]
and the initial prices $(t=0)$ are deduced from the distribution of $g_K$.
In the geometrical Brownian motion framework with $N_0=1$, the
Black--Scholes formula just computes the distribution of $g_1(N)$ as
$\p(g_1<t)={\cal N}(\sqrt t/2)-{\cal N}(-\sqrt t/2)=\p(4B_1^2\leq
t)$, where $B_1$ is a standard Gaussian random variable and ${\cal
N}(x)=\p(B_1\leq x)$ the Gaussian distribution function (see Profeta,
Roynette and Yor~\cite{PRY10} for a detailed study).

\subsubsection*{Financial framework} Assume $S$ to be a max-continuous nonnegative submartingale whose instantaneous return by time unit is
an adapted process $\lambda_t\geq0$ defined on a filtered probability spaced
$(\Omega,\mathcal F, (\mathcal{F}_t), \mathbb P)$. For instance,
$S$ is the current price of a stock under the risk neutral
probability in a financial market with short rate $\lambda_t$.
Put differently, $\tilde{S}_t=\exp(-\int_0^t\lambda_s\,ds)S_t$ is an
$(\mathcal{F}_t)$-martingale. We assume that $\int_0^\i\lambda_s\,ds=\i$ a.s.
Let $\zeta$ be an additional r.v. with conditional tail function
$\p(\zeta\geq t|\mathcal F_\i)=\exp(-\int_0^t\lambda_s \,ds)$.
Then $X_t=S_t\mathbf1_{t<\zeta}$ is a positive martingale with
negative jump to zero at time $\zeta$ with respect to the enlarged
filtration $\mathcal G_t=\sigma(\mathcal F_t,\zeta\wedge t).$
Since the $\mathbb G$-martingale $X$ goes to zero at $\infty$, the
random variable ${\overline X}_\zeta={\overline S}_\zeta$ is
distributed as $S_0/\xi$, where $\xi$ is uniformly distributed on
$[0,1]$. In particular, for any bounded function $h$
\[
\e[h({\overline S}_\zeta)]=
\e \biggl[\int_0^\infty e^{-\int_0^\alpha\lambda_sds}h({\overline S}_\alpha)\lambda_\alpha\, d\alpha \biggr]=
\int_0^1h(S_0/y)\,dy.
\]
In consequence we have access to the law of the properly discounted
maximum of the positive submartingale $S$. We stress that this is in
contrast to the more usual setting when one only has access to the
maximum of the discounted price process (cf. Grossman and Zhou~\cite{GZ93}).
We could also derive a conditional version of the equation above
representing $U(S_t)$ as a potential of the future maximum $\ovl
S_{t,u}$. Such representation find natural applications in financial
mathematics (see Bank and El Karoui~\cite{BankElKaroui04}).

\subsection{Az\'{e}ma--Yor martingales with given terminal values}
\label{sec:ay_given_values}
We describe now all local martingales whose terminal values are Borel functions
of the maximum of some nonnegative local martingale. This will be used
in subsequent sections, in particular to
construct Az\'{e}ma--Yor martingales with given terminal distribution
and solve the Skorokhod embedding problem.
We start with a simple lemma about solutions to a particular ODE.
%
\begin{lemma}
\label{lem:analytic} Let $h$ be a locally bounded Borel function
defined on $\re^+$, such that $h(x)/x^2$ is integrable away from zero.
Let $U$ be the solution of the ordinary differential equation (ODE),
%
%
\begin{equation}
\label{eq:odeayequation}
\forall x>0\qquad U(x)-xU'(x)=h(x)    \qquad \mbox{such that } \lim_{ x
\rightarrow
\i}U(x)/x=0 .
\end{equation}
%

 \textup{(a)}  The solution $U$ is given by
%
%
\begin{equation}
\label{eq:Ui}
U(x)=x \int_x^\i\frac{h(s)}{s^2}\,ds =\int_0^1 h\biggl(\frac{x}{s}\biggr)\,ds,\qquad x>0.
\end{equation}
%

 \textup{(b)}  Let $h_b(x) :=h(x\wedge b)$ be constant on $(b,\i)$.
The associated solution $U_b(x)=\int_0^1 h(\frac{x}{s}\wedge b)\,ds=U_b(x\wedge b)$
is constant on $(b,\i)$, and $U_b(x)=h_b(x)=h(b)$.

 \textup{(c)}  Let $h(m,x)=h(x\vee m)$ be constant on $(0,m)$. The
associated solution $U(m,x)$ is affine $U(m,x)=U(m)-U'(m)(m-x)$ for
$x\in(0,m)$.\vadjust{\goodbreak}
\end{lemma}

%
%
\begin{rem} \label{rem:moreonU}
We considered here $U$ on $(0,\infty)$ but naturally if $h$ is only
defined for $x>a>0$ then we consider $U$ also only for $x>a>0$. Note
that to define $U_b$ it suffices to have a locally integrable $h$
defined on $(0,b]$. We then put $h(x)=h(b)$, $x>b$.\vspace*{-3pt}
\end{rem}

\begin{pf*}{Proof of Lemma~\ref{lem:analytic}}
Formula \eqref{eq:Ui} is easy to obtain using the
transformation $(U(x)/x)'=-h(x)/x^2$ and the asymptotic condition in
\eqref{eq:odeayequation}.\break Both~(b) and (c) follow simply from \eqref{eq:Ui}.\vspace*{-3pt}
\end{pf*}

This analytical lemma allow us to characterize Az\'{e}ma--Yor
martingales from their terminal values. This extends, in more
detail, the ideas presented in El Karoui and Meziou
\cite{ELKM08}, Proposition~5.8.\vspace*{-3pt}
%
\begin{prop} \label{prop:terminalvalues}
Consider $(N_t)$ satisfying \eqref{eq:Ndef} with $N_0>0$ a \mbox{constant}.\vspace*{-2pt}

\begin{longlist}[(b)]
\item[(a)]
Let $h$ be a Borel function such that $h(x)/x^2$ is integrable away
from~$0$,
and $U$ the solution of the ODE \eqref{eq:odeayequation} given via
\eqref{eq:Ui}. Then
$h(\ovl N_\i)$ is an integrable random variable and the closed
martingale $\e (h(\ovl
N_\i)|\F_{t} )$, $t\geq0$, is the Az\'{e}ma--Yor martingale $M^{U}(N)$.
%
%
 \item[(b)] For a function $U$ with locally bounded derivative $U'$ and with
$U(x)/x\to0$ as $x\to\i$, the Az\'{e}ma--Yor local martingale
$M^U(N)$ is a uniformly integrable martingale if and only if
$h(x)/x^2$ is integrable away from zero,
where~$h$ is now defined via \eqref{eq:odeayequation}.\vspace*{-3pt}
\end{longlist}
\end{prop}
\begin{pf}
We start with the proof of (a). We have
\[
\e|h(\ovl N_\i)|=\int_0^1 |h(N_0/s)|\,ds =N_0 \int_{N_0}^\i
|h(s)|/s^2\, ds<\i,
\]
since we assumed that $h(x)/x^2$ is integrable away from $0$.
To study the martingale $H_t=\e (h(\ovl N_\infty)|\F_t)$, we use
that $\ovl N_\i= \ovl N_t \vee\ovl N_{t,\infty}$. From
Proposition~\ref{prop:uniform}, the running maximum $\ovl N_{t,\i}$,
conditionally on $\F_t$, is distributed as $N_t/\xi$, for an
independent r.v. $\xi$ uniform on $[0,1]$.
The martingale $H_t$ is given by the following closed formula $H_t=\e
 (h (\ovl N_t\vee(N_t/\xi))|\F_t )$, that is,
\[
H_t=\int_0^1h\bigl(\ovl N_t\vee(N_t/s)\bigr)\,ds=U(\ovl N_t,N_t)=U(\ovl
N_t)-U'(\ovl N_t)(\ovl N_t-N_t),
\]
where in the
last equalities, we have used  (c)  in Lemma~\ref{lem:analytic}.

Given (a), to prove part (b) it suffices to observe that $M_t^U(N)\to
h(\ovl N_\i)$ a.s.
and hence integrability of $h(\ovl N_\i)$, that is, integrability of
$h(x)/x^2$
away from zero, is necessary for uniform integrability of $M^U(N)$.\vspace*{-3pt}
\end{pf}

%
\begin{rem}
It is not necessary to assume that $N_0$ is a constant in Proposition
\ref{prop:terminalvalues}.
However, if $N_0$ is random we have to further assume that $\e\int
_0^1 |h(N_0/s)|\,ds=\int_1^\i\e|h(xN_0)|\frac{dx}{x^2}<\i$.
This holds, for example, if $N_0$ is integrable and $N_0>\epsilon>0$
a.s. We can apply the same reasoning to the process $(N_{t+u}\dvtx  u\geq
0)$ to see that if $\e\int_0^1 |h(N_t/s)|\,ds<\i$ then $U(N_t)=\e
 (h(\ovl N_{t,\i})|\F_{t} )$.\vadjust{\goodbreak}

Finally, we note that similar consideration as in  (a) above were
independently made in Nikeghbali and Yor~\cite{NY06}.

We stress that the boundary condition $U(x)/x\to0$ as $x\to\i$ for
\eqref{eq:odeayequation} is essential in part  (a).
Indeed, consider $N_t=1/Z_t$ the inverse of a three-dimensional Bessel process.
Note that $N_t$ satisfies our hypothesis, and it is well known that $N_t$
is a strict local martingale (cf. Exercise V.2.13 in Revuz and Yor
\cite{ry}).
Then for $U(x)=x$ we have $M^U_t(N)=N_t$ is also a strict local martingale,
but obviously we have $U(x)-U'(x)x=0$.

Observe that $\p(N_\zeta\in\{0,b\})=1$ if $\zeta=\inf\{t\geq0\dvtx
N_t\notin(0,b)\}$. Then, if~$h$ is constant on
$[b,\i)$, then $h(\ovl N_\i)=h(\ovl N_{\zeta})$, and the closed
martingale $\e (h(\ovl N_{\zeta})|\break\F_{t\wedge\zeta} )$,
$t\geq0$, is the Az\'{e}ma--Yor martingale
$M^{U_b}_{t}(N)=M^{U_b}_{t\wedge\zeta}(N)$, where $U_b$ the solution
of the ODE \eqref{eq:odeayequation}
given in point (b) of Lemma~\ref{lem:analytic}.\vspace*{-3pt}
\end{rem}

As shown in Sections~\ref{sec:def} and~\ref{sec:bach}, Az\'{e}ma--Yor
processes $Y=M^U(N)$ generated by an increasing function $U$ have very
nice properties based on the characterization of their maximum as $\ovl
Y=U(\ovl N)$. In particular, from Theorem~\ref{thm:DD}, the process
$Y$ satisfies a DD-constraint and can also be characterized from its
terminal value. Recall Definitions~\ref{def:DD},~\ref{def:w-DD} and
the stopping time $\zeta_w^Y$ defined before Theorem~\ref{thm:DD}.\vspace*{-3pt}
%
\begin{prop}
\label{prop:increasing_function}
Let $h$ be a right-continuous nondecreasing function such that
$h(x)/x^2$ is integrable away from $0$, and put $b=\inf\{x\dvtx
h(y)=h(x)\ \forall y\geq x\}$.\vspace*{-2pt}
\begin{longlist}[(b)]
\item[(a)] The solution $U$ of the ODE \eqref{eq:odeayequation} is then a
continuous strictly increasing concave function on $(0,b)$ and constant
and equal to $h(b)$ on $(b,\i)$.
\item[(b)] Let $V$ be the inverse of $U$ with $V(U(b))=b$. Function
$w(y)=h(V(y))$ given by \eqref{eq:DDfunction} for $y<U(b)$ and by
$w(y)=y$ for $y\geq U(b)$, is a right-continuous drawdown function and
$r_w=U(b)=h(b)$.
\item[(c)] Consider $(N_t)$ satisfying \eqref{eq:Ndef} with $N_0>0$ a
constant. The uniformly integrable martingale $Y_t=M^U_t(N)=\e[h(\ovl
N_\i)|\F_t]$ satisfies $w$-DD constraint. Furthermore, $Y_t=Y_{t\land
\zeta_w^Y}$, $Y_{\zeta_w^Y}=w(\ovl Y_{\zeta_w^Y})$ a.s. and $\zeta
_w^Y=\inf\{t\dvtx  N_t\notin(0,b)\}$.\vspace*{-2pt}
\end{longlist}

Conversely, let $w$ be a right-continuous drawdown function, with
functions $V,U,h$ satisfying \textup{(a)} and \textup{(b)}. Then any uniformly
integrable
martingale $Y$, satisfying the $w$-DD constraint and $Y_{\zeta
_w^Y}=w(\ovl Y_{\zeta_w^Y})$ a.s., is an Az\'{e}ma--Yor martingale
$M^U(N)$ for some $(N_t)$ satisfying \eqref{eq:Ndef} with
$N_0=V(Y_0)>0$ and such that $N_{t\land\zeta_w^Y}=M^V_{t\land\zeta
_w^Y}(Y)$ and $\zeta_w^Y=\inf\{t\dvtx  N_t\notin(0,V(r_w)\}$.\vspace*{-3pt}
\end{prop}

%
\begin{rem}
Note that $h$, and consequently $U$, need to be defined only for $x\geq
N_0$. Then $V(y)$ is defined for $U(N_0)\leq y\leq U(b)$ with
$V(U(b))=b\leq\i$ and the drawdown function $w(y)$ is defined for
$y\geq U(N_0)$.

A solution $U$ of the ODE \eqref{eq:odeayequation} is strictly
increasing if and only if
$U>h$; however, only increasing and concave solutions are easy to
characterize.\vadjust{\goodbreak}
\end{rem}

\begin{pf*}{Proof of Proposition~\ref{prop:increasing_function}}
(a) When $h$ is nondecreasing, from \eqref{eq:odeayequation} and
\eqref{eq:Ui} it is
clear that $U$ is strictly increasing until that $h$ becomes constant,
and constant after that. If $h$ is differentiable, concavity of $U$
follows since $-xU''(x)=h'(x).$ The general case follows by
regularization or can be checked directly using \eqref{eq:Ui} which
yields to $U'(x)=\int_x^\i(h(s)-h(x))/s^2\,ds=\int_0^\i
(h(s)-h(x))^+/s^2\,ds$.

(b) In consequence, $V$ is increasing and convex on $[U(0),U(b))$, and
hence by \eqref{eq:DDfunction} $w(y)$ is increasing and $w(y)<y$ for
$y\in(U(0),U(b))$. We thus have $r_w=U(b)$ but note that we could have
$w(U(b)-)$ both less then or equal to $U(b)$. Integrability properties
of $w$ in Definition~\ref{def:w-DD} follow since $V$ and $U$ are well
defined and we conclude that $w$ is a drawdown function.
Right-continuity of $w$ follows from right-continuity of $h$.
More precisely, from~\eqref{eq:odeayequation}, $u=U'$ is right
continuous, and hence also $V'(y)=1/u(V(y))$ is right-continuous and
nondecreasing.

(c) Identification of $Y$ is given in part (a) of Proposition \ref
{prop:terminalvalues}.
The rest follows from Lemma~\ref{lem:Nstop}, Theorem~\ref{thm:DD} and
the fact
that $N_\zeta\in\{0,b\}$ a.s. for $\zeta=\inf\{t\dvtx  N_t\notin(0,b)\}
$ upon noting that $V(r_w)=b$.
\end{pf*}


\subsection{\texorpdfstring{On relations between $\operatorname{AVaR}_\mu$, Hardy--Littlewood transform and tail quantiles}
{On relations between AVaR mu, Hardy--Littlewood transform and tail quantiles}}
\label{sec:AVaR_HL}
In this section we present results about probability measures, their
tail quantile function, the average value at risk and the
Hardy--Little\-wood transform.
The presentation is greatly simplified using tail quantiles of measure.

The notation and quantities now introduced will be used throughout the
rest of the paper. For $\mu$ a probability measure on $\re$ we
denote $l_\mu,r_\mu$, respectively, the lower and upper bound of the
support of $\mu$.
We let $\ovl\mu(x)=\mu ([x,\i) )$ and
$\ovl q_\mu\dvtx  (0,1]\to\re\cup\{\i\}$ be the tail quantile function
defined as the left-continuous inverse of $\ovl\mu$, $\ovl
q_\mu(\lambda):= \inf\{x\in\re\dvtx   \ovl\mu(x) < \lambda\}$.
When ${\ovl q}_\mu(\lambda)$ is a point of
continuity of $\ovl\mu$, then ${\ovl\mu}({\ovl q}_\mu(\lambda
))=\lambda$, whereas if not, ${\ovl\mu}({\ovl q}_\mu(\lambda
)+)<\lambda\leq{\ovl\mu}({\ovl q}_\mu(\lambda))$. In particular,
if ${\ovl\mu}(r_\mu)>0$, $r_\mu=\ovl q(0+)$ is a
jump of ${\ovl\mu}$ and $r_\mu=\ovl q(0+)=\ovl q(\lambda)$, if
$0<\lambda\leq{\ovl\mu}(r_\mu)$.

We write $X\sim\mu$ to denote that $X$ has distribution
$\mu$ and recall that ${\ovl q}_\mu(\xi)\sim\mu$ for $\xi$ uniformly
distributed on $[0,1]$.

\setcounter{footnote}{2}

Assume $\int_{\re} |s|\mu(ds)<\i$ and let
$m_\mu=\int_{\re} s\mu(ds)$. We define call
function\footnote{This denomination is used
in financial literature while the actuarial literature uses rather the
notion of
stop-loss function (cf.~\cite{Hu98}).}~$C_\mu$ and barycenter
function $\psi_\mu$ by
%
%
\begin{equation}\label{eq:call+barycentre}
{C}_\mu(K)=\int_{\re}(s-K)^+\mu(ds),\qquad
\psi_{\mu}(x)=
\frac{1}{{\ovl\mu}(x)}\int_{[x,\infty)}s\mu(ds),
\end{equation}
where $K\in\re,x<r_\mu$. We put $\psi_\mu(x)=x$ for $x\geq 
r_\mu$.\vadjust{\goodbreak}

Finally, we also introduce the \textit{average value at risk} at the
level $\lambda\in(0,1]$,
given by
%
%
\begin{equation}\label{eq:AVaR}
\operatorname{AVaR}_\mu(\lambda)=\frac{1}{\lambda}\int_0^\lambda\ovl q_\mu
(\eta)\,d\eta .
\end{equation}
Observe that $\operatorname{AVaR}_\mu$ is equal to $r_\mu={\ovl q}_\mu(0+)$
on $(0,\ovl\mu(r_\mu)]$,
is strictly decreasing on $({\ovl\mu}(r_\mu),1)$ since its
derivative is then equal to $\frac{1}{\lambda^2}\int_0^\lambda(\ovl
q_\mu(\lambda)-\ovl q_\mu(\eta))\,d\eta<0$,
and $\operatorname{AVaR}_\mu(1)=m_\mu$.

The average value at risk $\operatorname{AVaR}_\mu$ is thus a quantile
function of some probability measure $\mu^{\mathrm{HL}}$ with support
$(m_\mu,r_\mu)$, which can be defined by
%
%
\begin{equation}
\label{eq:HL_definition} \mu^{\mathrm{HL}}\sim\operatorname{AVaR}_\mu(\xi),\qquad
\xi\mbox{ uniform on }[0,1].
\end{equation}
This distribution, called the \textit{Hardy--Littlewood transform} of
$\mu$ has been intensively studied by many authors, starting with the
famous paper of Hardy and Littlewood~\cite{HL30}. We will describe its
prominent role in the study of distributions of maxima of martingales
in Section~\ref{sec:optimal_prop} below.
Recently F\"{o}llmer and Schied~\cite{FS04} studied properties of
$\operatorname{AVaR}_{\mu}$ as a coherent risk measure. Finally note that here
$\mu$ is the law of losses (i.e., negative of gains) and some
authors refer to $\operatorname{AVaR}_\mu(\lambda)$
as $\operatorname{AVaR}_\mu(1-\lambda)$.

As in~\cite{FS04}, pages 179--182, page 408, Lemma A.22, it is easy to characterize the Fenchel transform of
the concave function $\lambda\operatorname{AVaR}_\mu(\lambda)$ as the
call function. From this property, we infer a nonclassical
representation of the tail function $\ovl\mu^{\mathrm{HL}}(y)$ as an
infimum.
%
\begin{prop} \label{prop:AVaR}
Let $\mu$ be a probability measure on $\re$, \mbox{$\int|s|\mu(ds)<\i
$}.

\begin{longlist}[(b)]
\item[(a)] The average value at risk $\operatorname{AVaR}_\mu(\lambda)$ can be
described as, $\lambda\in(0,1)$,
%
%
\begin{equation}
\label{eq:AVaR_2}
\operatorname{AVaR}_\mu(\lambda) =\frac{1}{\lambda}C_\mu(\ovl
q_\mu(\lambda))+\ovl q_\mu(\lambda)=\frac{1}{\lambda}\inf_{K\in
\re}\bigl(C_\mu(K)+\lambda
K\bigr).
\end{equation}
\item[(b)] The call function is the Fenchel transform of
$\lambda\operatorname{AVaR}_\mu(\lambda)$, so that
%
%
\begin{equation}
\label{eq:calls_via_AVaR}
C_\mu(K)=\sup_{\lambda\in(0,1)} \bigl(\lambda\operatorname{AVaR}_\mu(\lambda
)-\lambda K\bigr),\qquad K\in\re.
\end{equation}
\item[(c)] The Hardy--Littlewood tail function ${\ovl\mu}^{\mathrm{HL}}$ is given
for any $y\in(m_\mu,r_\mu) $ by
%
%
\begin{equation}
\label{eq:HL_dual}
{\ovl\mu}^{\mathrm{HL}}(y)=\inf_{z>0} \frac{1}{z}C_\mu(y-z).
\end{equation}
\item[(d)] The barycenter function and its right-continuous inverse are related
to the average value at risk and Hardy--Littlewood tail function by
%
%
\begin{eqnarray}
\label{eq:barycentre and dual}
\psi_{\mu}(x)&=&\operatorname{AVaR}_\mu({\ovl\mu}(x)), \qquad  x\leq
r_\mu,\nonumber
\\[-8pt]
\\[-8pt]
\psi_{\mu}^{-1}(y)&=&\ovl q_\mu( {\ovl\mu}^{\mathrm{HL}}(y)), \qquad  y\in[m_\mu,
r_\mu].
\nonumber
\end{eqnarray}
\end{longlist}
\end{prop}

%
\begin{rem}
From \eqref{eq:call+barycentre} we have $\psi_\mu(x)=\e[X|X\geq x]$,
where $X\sim\mu$. Then \eqref{eq:barycentre and dual} gives
$\operatorname{AVaR}_\mu(\lambda)=\e[X|X\geq\ovl q_\mu(\lambda)]$, $d\ovl q_\mu(\lambda)$-a.e.,
which justifies names \textit{expected shortfall}, or \textit{conditional
value at risk}
used for $\operatorname{AVaR}_\mu$.
\end{rem}

\begin{pf*}{Proof of Proposition~\ref{prop:AVaR}} We write $\ovl q=\ovl q_\mu$.
\begin{longlist}[(b)]
\item[(a)] The proof is based on the classical property, ${\ovl
q}(\xi)\sim\mu$, for $\xi$ uniformly distributed on $[0,1]$.
Then
\[
C_\mu(\ovl q(\lambda))=\int_0^1\bigl(\ovl q(\eta)-\ovl
q(\lambda)\bigr)^+\,
d\eta=
\int_0^\lambda\bigl(\ovl q(\eta)-\ovl q(\lambda)\bigr)\,d\eta= \lambda
\bigl(\operatorname{AVaR}_\mu(\lambda)-\ovl q(\lambda)\bigr).
\]
Moreover, the convex function $G_\lambda(K):=C_\mu(K)+\lambda K$
attains its minimum in $K_\lambda$ such that
$\ovl\mu(K_\lambda)=\lambda$.

When ${\ovl\mu}({\ovl q}(\lambda))=\lambda$, ${\ovl q}(\lambda)$
is a minimum of the function $G_\lambda(K)$, $\lambda\operatorname{AVaR}_\mu(\lambda)=G_\lambda(\ovl q(\lambda))$, and
\eqref{eq:AVaR_2} holds true.

If ${\ovl\mu}({\ovl q}(\lambda))>\lambda> {\ovl\mu}({\ovl
q}(\lambda)+)$, then $\mu$ has an atom in $x:={\ovl q}(\lambda)$.
$G_\lambda$ has a~minimum in $x$ and
$G'_\lambda$ changes sign discontinuously in $x$. Then we see
that $G_\lambda(\ovl q(\lambda)=G_\lambda(x)$ is linear in
$\lambda\in(\ovl\mu(x+),\ovl\mu(x))$.
\item[(b)] Convex duality for Fenchel transforms yields
\eqref{eq:calls_via_AVaR} from \eqref{eq:AVaR_2}.
\item[(c)]
Using \eqref{eq:AVaR_2} we have, for any $y>m_\mu$ and $\lambda\in(0,1)$
%
%
\begin{eqnarray}
\label{eq:HLCaracterisation}
 \operatorname{AVaR}(\lambda)<y    \quad &\Leftrightarrow& \quad \exists K     \mbox{ such
that }  y>\frac{C_\mu(K)}{\lambda}+K\nonumber
\\
   &\Leftrightarrow& \quad \exists K<y   \mbox{ such that }
\lambda>\frac{C_\mu(K)}{y-K}  \\
 &\Leftrightarrow& \quad
\lambda>\inf_{K<y}\frac{C_\mu(K)}{y-K}.
\nonumber
\end{eqnarray}
The function
$\inf_{K<y}\frac{C_\mu(K)}{y-K}=\inf_{z>0}\frac{1}{z}C_\mu(y-z)
$ is decreasing and left continuous. We conclude that it is the
left-continuous inverse function of $\operatorname{AVaR}_\mu(\lambda)$ which
is $\ovl\mu^{\mathrm{HL}}$.
\item[(d)] By definition, $\ovl\mu(x)\operatorname{AVaR}_\mu(\ovl\mu(x))=\int
_0^{\ovl\mu(x)} \ovl q_\mu(\eta)\,d\eta=\int_{[x,\infty)}s\mu(ds)$.

The right-continuous inverse $\psi_{\mu}^{-1}(y)$ of the nondecreasing left-continuous function
$\psi_{\mu}$ is defined by $\psi_{\mu}^{-1}(y)=\sup\{x\dvtx \psi_{\mu
}(x)\leq y\}=
\sup\{x\dvtx\break \operatorname{AVaR}_\mu({\ovl\mu}(x))\leq y\}$. Since, $\ovl\mu
^{\mathrm{HL}}$ is the left-continuous inverse of
$\operatorname{AVaR}_\mu$, the following inequalities hold true for $y\in
[m_\mu,r_\mu]\dvtx\psi_{\mu}^{-1}(y)=\sup\{x\dvtx {\ovl\mu}(x)
\geq\break\ovl\mu^{\mathrm{HL}}(y)\}=\sup\{x\dvtx  x\leq\ovl q(\ovl\mu^{\mathrm{HL}}(y))\}
=\ovl q(\ovl\mu^{\mathrm{HL}}(y)).$
\qed\end{longlist}
\noqed\end{pf*}

We now describe the relationship between $\mu$, $\operatorname{AVaR}_\mu$,
$\psi_\mu$ and $\mu^{\mathrm{HL}}$ on one hand, and $w_\mu$, solutions
$U_\mu$ of \eqref{eq:odeayequation} when $h(x)=\ovl q_\mu(1/x)$ and
the associated Az\'{e}ma--Yor martingales $M^{U_\mu}(N)$ on the other
hand. It turns out all these objects are intimately linked together in
a rather elegant manner. Some of our descriptions below, in particular
characterization of $\operatorname{AVaR}$ in (a), appear to be different from
classical forms in the literature. We note that we start with $\mu$
and define $h$ but equivalently we could start with a nondecresing
right-continuous $h$ and use $h(x)=\ovl q_\mu(1/x)$ to define $\mu$.

Recall Definitions~\ref{def:DD},~\ref{def:w-DD} and the stopping time
$\zeta_w^Y$ defined before Theorem~\ref{thm:DD}.
%
\begin{prop}
\label{prop:q_to_U} Let $\mu$ be a probability measure on $\re$,
\mbox{$\int|s|\mu(ds)<\i$}.

\begin{longlist}[(b)]
\item[(a)] $U_\mu(x):=\operatorname{AVaR}_\mu(1/x)$ solves \eqref{eq:odeayequation}
with $h_\mu(x)=\ovl q_\mu(1/x)$, $x\geq1$, and $h_\mu(x)/x^2$ is
integrable away from zero. In particular $U_\mu$ is given by \eqref
{eq:Ui} and $U_\mu(x)=U_\mu(x\wedge b_\mu)$ with $b_\mu=1/\ovl\mu
(r_\mu)$. $U_\mu$
is concave and $V_\mu(y)=1/\ovl\mu^{\mathrm{HL}}(y)$ is the inverse function
of $U_\mu$.
\item[(b)] Let $w_\mu(y)=h_\mu(V_\mu(y))=q_\mu(\ovl\mu^{\mathrm{HL}}(y))$ be the
function associated with~$\mu$ by \eqref{eq:DDfunction} for $y\in
(m_\mu,r_\mu)$, and extended via $w_\mu(y)=y$ for $y\geq r_\mu$.
Then $w_\mu$ is a drawdown function, $r_w=r_\mu$ and $w_\mu$ is the
right-continuous inverse of the barycenter function~$\psi_\mu$.
Furthermore, $w_\mu$ is the hyperbolic derivative of~$V_\mu$ as
defined by Kertz and R\"{o}sler~\cite{KertzRosler93}.
\item[(c)] Let $N$ satisfy \eqref{eq:Ndef} with $N_0=1$ and $Y_t=M_{t}^{U_\mu
}(N)$. Then $Y_t\geq U_\mu(N_t)$, $Y_\i=Y_{\zeta_{w_\mu}^Y}=\ovl
q_\mu (1/\ovl N_\zeta )$ is distributed according to $\mu$
and $\ovl Y_\infty=U_\mu(\ovl N_\zeta)=\operatorname{AVaR}_\mu(1/\ovl
N_\zeta)$ is distributed according to $\mu^{\mathrm{HL}}$. Furthermore, the
process $(Y_t)$ is a uniformly integrable martingale which satisfies
$w_\mu$-DD constraint and $\zeta_{w_\mu}^Y=\inf\{t\dvtx  N_t\notin
(0,b_\mu)\}=\inf\{t\dvtx  \psi_\mu(Y_t)\leq\ovl Y_t\}$.
\end{longlist}

The same properties hold true for any max-continuous uniformly
integrable martingale $Y$, $Y_0=U(1)$, satisfying the $w_\mu$-DD
constraint up to $\zeta=\zeta^Y_{w_\mu}$ and $Y_{\zeta}=w_\mu(\ovl
Y_{\zeta})$ a.s.
\end{prop}

\begin{pf} We write $\ovl q=\ovl q_\mu$, $r=r_\mu$, $b=b_\mu
=1/\ovl\mu(r_\mu)=V_\mu(r_\mu-)$.
{\smallskipamount=0pt\begin{longlist}[(b)]
\item[(a)] From definition \eqref{eq:AVaR} we have $\operatorname{AVaR}_\mu(\lambda
)=\int_0^1 \ovl q(\lambda s)\,ds $ which is exactly formula \eqref
{eq:Ui}. Note that in the case $\ovl\mu(r)>0$ we have $h(x)=h(b)$,
$x\geq b$ with $b=1/\ovl\mu(r)$. $U_\mu$ is concave by Proposition
\ref{prop:increasing_function}. The rest follows since $\operatorname{AVaR}_\mu
$ is the tail quantile of $\mu^{\mathrm{HL}}$ [see \eqref{eq:HL_definition}].
\item[(b)] This follows by part (d) in Proposition~\ref{prop:AVaR} and the last
statement follows from \eqref{eq:DDEquation} and Theorem~4.3 in \cite
{KertzRosler93}.
\item[(c)] We have $Y_t\geq U(N_t)$ from concavity of $U_\mu$. The rest
follows easily from points (a) and (b) above together with Proposition
\ref{prop:increasing_function}, properties of $\ovl q$, universal law
of $\ovl N_\zeta$ given in point (c) of Proposition~\ref{prop:uniform}
and the definition of $\mu^{\mathrm{HL}}$ in \eqref{eq:HL_definition}.
\qed
 \end{longlist}}
\noqed\end{pf}

\subsubsection*{An illustrative example (continued from Section \protect\ref
{sec:dd})}
We come back to the example with linear DD-constraint
$w(y)=\gamma y$, $0<\gamma<1$, resulting from function $U(x)=\frac
{1}{1-\gamma} x^{1-\gamma}$, $x\geq1$. Using Proposition \ref
{prop:q_to_U} we have $Y_\i\sim\mu$ and $\ovl Y_\i\sim\mu^{\mathrm{HL}}$
which we can now easily describe. We have $\ovl\mu
^{\mathrm{HL}}(y)=1/V(y)=((1-\gamma)y)^{1/(\gamma-1)}$ for $y\geq m_\mu=\operatorname{AVaR}_\mu(1)=U(1)=\frac{1}{1-\gamma}$.
In consequence, the random
variable $\ovl Y_\i$ is distributed according to
a \textit{Pareto distribution}, with shape parameter $a=m_\mu=\frac
{1}{1-\gamma}$ and location parameter $m=m_\mu$. The mean of $\ovl
Y_\i$ is $a m/(a-1)=(\gamma(1-\gamma))^{-1}$. Since $Y_\i=w(\ovl
Y_\i)=\gamma\ovl Y_\i$ we see that $Y_\i$ is still distributed
according to a Pareto distribution, with the same shape parameter, and
location parameter $m_{1}=m \gamma=\frac{\gamma}{1-\gamma}$.
Naturally, we could also describe $\mu$ using $\ovl q_\mu(\lambda
)=h(1/\lambda)=\frac{\gamma}{1-\gamma}\lambda^{\gamma-1}$ which,
taking inverses, gives $\ovl\mu(x)= (\frac{\gamma}{1-\gamma
}\frac{1}{x} )^{\fracc{1}{1-\gamma}}$ as required.

As a consequence we see that if $(Y_t)$ is a max-continuous martingale
which satisfies a linear constraint $Y_t> \gamma\ovl Y_t$ until $\zeta
=\zeta_w^Y<\i$ a.s., then, necessarily, $Y_{\zeta}=\gamma\ovl
Y_{\zeta}$ has a Pareto distribution.

\subsection{The Skorokhod embedding problem revisited}
\label{sec:skoro}
The Skorokhod embedding problem can be phrased as follows:
given a probability measure $\mu$ on~$\re$
find a stopping time $T$ such that $X_T$ has the law $\mu$,
$X_T\sim\mu$. One further requires $T$ to be \textit{small}
in some sense, typically saying that $T$ is minimal.
We refer the reader to Ob\l\'{o}j~\cite{genealogia}
for further details and the history of the problem.

In~\cite{ay1} Az\'{e}ma and Yor introduced the family of martingales
described in Definition~\ref{def:ay} and used them to give an
elegant solution to the Skorokhod embedding problem for $X$ a
continuous local martingale (and $\mu$ centered). Namely, they
proved that
%
%
\begin{equation}\label{eq:aystop}
T_{\psi}(X)=\inf\{t\geq0\dvtx  \psi_\mu(X_t)\le\ovl X_t\},
\end{equation}
where $\psi_\mu$ is the barycenter function \eqref
{eq:call+barycentre}, solves the embedding problem.

We propose to rediscover their solution in a natural way using our
methods, based on the observation that the process $X$ satisfies the
$w_\mu$-DD constraint up to $T_{\psi}(X)$. If we show the equality
$X_\zeta=w_\mu(\ovl X_\zeta)$ at time $\zeta=T_{\psi}(X)$,
Proposition~\ref{prop:q_to_U} gives us the result.
%
\begin{theorem}[(Az\'{e}ma and Yor~\cite{ay1})]
\label{th:AYor} Let $(X_t)$ be a continuous local martingale, $X_0\in
\re$ a constant, $\langle X\rangle_\i=\i$ a.s. and $\mu$ a
probability measure on $\re\dvtx\int|x|\mu(dx)<\i$, $\int x\mu
(dx)=X_0$. Then $T_\psi<\i$ a.s., $(X_{t\land T_{\psi}})$ is a
uniformly integrable martingale and $X_{T_{\psi}}\sim\mu,  \ovl
X_{T_{\psi}}\sim\mu^{\mathrm{HL}}$, where $T_{\psi}$ is defined via \eqref
{eq:aystop}.

With notation of Proposition~\ref{prop:q_to_U}, define $N_t=M^{V_\mu
}_{t\land\tau^{r_\mu}(X)}(X)$. Then
%
%
\begin{equation}\label{eq:AYisDD}
T_{\psi} = \inf\{t\geq0\dvtx  X_t\leq w_\mu(\ovl X_t)\}=\inf\{t\geq0\dvtx
N_t\leq0\}\wedge\tau^{b_\mu}(N)
\end{equation}
and $X_{t\wedge T_\psi}=M^{U_\mu}_{t\wedge T_\psi}(N)$.
\end{theorem}

\begin{pf} Let $\tau=\tau^{r_\mu}(X)$. $(N_t\dvtx t<\tau)$ is a
continuous local martingale with $N_0=1$ since $U_\mu(1)=\operatorname{AVaR}_\mu(1)=X_0$.
If $b_\mu<\i$, then $r_\mu<\infty$ and $(N_t\dvtx t\leq\tau)$ is a
local martingale stopped at $\inf\{t\dvtx  N_t=b_\mu\}=\tau<\i$ a.s.
Suppose $b_\mu=\i$. Then $\ovl N_{\tau-}=\lim_{x\to r_\mu}V(x)=\i
$. This readily
implies that $\langle N\rangle_{\tau-}=\i$ a.s. and in particular
$\tau_0(N)<\tau$ a.s.  (cf. Proposition~V.1.8 in Revuz and Yor \cite
{ry}). Note that this applies both for the case $r_\mu$ finite and
infinite. We conclude that $N_{t\land\tau_0(N)}$ is a continuous
local martingale satisfying \eqref{eq:Ndef} stopped at $\tau
_0(N)\wedge\tau^{b_\mu}(N)<\i$ a.s. The theorem now follows from
part (c) in Proposition~\ref{prop:q_to_U}.
\end{pf}

%
\begin{rem}\label{rem:skoro_maxcont}
Note that in general only max-continuity of $(X_t)$ would not be
enough. More precisely we need to have $X_{T_\psi}=w_\mu(\ovl
X_{T_\psi})$ a.s. or equivalently that the process $N_t$ crosses zero
continuously. Note also that we do \textit{not} necessarily have that
$\psi_\mu(X_{T_\psi})=\ovl X_{T_\psi}$. Finally, we point out that
the value of $X_0=\int x\mu(dx)$ plays no special role and we do
not need to assume that $X_0=0$.
\end{rem}

\section{On optimal properties of AY martingales related to HL
transform and its inverse}
\label{sec:optimal_prop}
In this final section we investigate the optimal properties of Az\'
{e}ma--Yor processes and of the Hardy--Littlewood transform $\mu\to\mu
^{\mathrm{HL}}$ and its inverse operator $\Delta$. We use two orderings of
probability measures. We say that $\mu$ dominates $\nu$ in the \textit
{stochastic order} (or stochastically) if $\ovl\mu(y)\geq\ovl\nu
(y)$, $y\in\re$. We say that $\mu$ dominates $\nu$ in the \textit
{increasing convex order} if $\int g(y)\mu(dy)\geq\int g(y)\nu
(dy)$ for any increasing convex function $g$ whenever the integrals
are defined.
Observe that the latter order is equivalent to $C_\mu(K)\geq C_\nu
(K)$, $K\in\re$
(cf. Shaked and Shanthikumar~\cite{ShakedShanthikumar94}, Theorem~3.A.1).

From \eqref{eq:AVaR_2} and \eqref{eq:calls_via_AVaR} we deduce
instantly that if $\mu,\rho$ are probability measures on $\re$ which
admit first moments, then
%
%
\begin{eqnarray} \label{eq:AVaR_to_calls}
&&C_\mu(K)\leq C_\rho(K), \qquad  K\in\re\nonumber
\\[-8pt]
\\[-8pt]
&&\qquad \Leftrightarrow \quad \operatorname{AVaR}_\mu
(\lambda)\leq\operatorname{AVaR}_\rho(\lambda), \qquad  \lambda\in(0,1).
\nonumber
\end{eqnarray}
By definition of $\mu^{\mathrm{HL}}$, $\operatorname{AVaR}_\mu(\lambda)=\ovl q_{\mu
^{\mathrm{HL}}}(\lambda)$, and hence we obtain
%
%
\begin{eqnarray}
\label{eq:HL_to_calls}  \hspace*{18pt}
C_\mu(K)\leq C_\rho(K), \qquad  K\in\re &  \quad  \Leftrightarrow \quad &
{\ovl q}_{\mu^{\mathrm{HL}}}(\lambda)\leq\ovl q_{\rho^{\mathrm{HL}}}(\lambda), \qquad
\lambda\in[0,1] \nonumber
\\[-8pt]
\\[-8pt]
\hspace*{18pt} & \Leftrightarrow& {\ovl\mu}^{\mathrm{HL}}(y)\leq{\ovl\rho}^{\mathrm{HL}}(y), \qquad
y\in\re,
\nonumber
\end{eqnarray}
so that $\rho^{\mathrm{HL}}$ dominates $\mu^{\mathrm{HL}}$ stochastically if and only
if $\rho$ dominates $\mu$ in the convex order.
\subsection{Optimality of Az\'{e}ma--Yor stopping time and
Hardy--Littlewood transformation}

The Az\'{e}ma--Yor stopping time has
a remarkable property that the maximum of a martingale stopped at this
time is distributed according to $\mu^{\mathrm{HL}}$ [see Theorem~\ref
{th:AYor}, also part  (c) in Proposition~\ref{prop:q_to_U}]. The
importance of this result comes from the result of Blackwell and Dubins
\cite{BD63} (see also the concise version of Gilat and Meljison \cite
{GM88}) showing:
%
\begin{theorem}[(Blackwell and Dubins~\cite{BD63})]\label{thm:BD} Let $(P_t)$ be
a uniformly integrable martingale and $\mu$ the distribution of
$P_\infty$. Then
%
%
\begin{equation}
\label{eq:BD_domination}
\p(\ovl P_\infty\geq y)\leq\mu^{\mathrm{HL}}([y,\infty)),\qquad y\in\re.
\end{equation}
In other words, any Hardy--Littlewood maximal r.v. associated with
$P_\infty$ dominates stochastically $\ovl P_\infty$.
\end{theorem}

 In fact $\mu^{\mathrm{HL}}$ is sometimes defined as the smallest
measure which satisfies~\eqref{eq:BD_domination}. One then proves the
representation \eqref{eq:HL_definition}.
\begin{pf*}{Proof of Theorem~\ref{thm:BD}}
We present a short proof of the theorem based on arguments in Brown,
Hobson and Rogers~\cite{BHR01}.
Let $(P_t)$ be a uniformly integrable martingale with terminal distribution $\mu$. Choose $y\in (0,r_\mu)$ and observe that for any $K<y$ the
following inequality holds a.s.:
%
%
\begin{equation}
\label{eq:bhr}
\mathbf{1}_{\ovl P_\infty\geq y}\leq\frac{(P_\infty
-K)^+}{y-K}+\frac{y-P_\infty}{y-K}\mathbf{1}_{\ovl P_\infty\geq y} .
\end{equation}
If $P$ is max-continuous then the last term on the RHS is simply
$-M^F_\infty(P)$ for $F(z)=\frac{(z-y)^+}{y-K}$ and has zero
expectation. In general, we can substitute the last term with a greater
term $\frac{P_{\tau^y(P)}-P_\infty}{y-K}\mathbf{1}_{\ovl P_\infty
\geq y}$ which has zero expectation.
Hence, taking expectations in \eqref{eq:bhr} we find
\[
\p(\overline{P}_\i\geq y)\leq\frac{1}{y-K}\int_K^\infty(x-K)\mu
(dx).
\]
Taking infimum in $K<y$ and using \eqref{eq:HL_dual} we conclude that
$\p(\overline{P}_\i\geq y)\leq\ovl\mu^{\mathrm{HL}}(y)$.
\end{pf*}

Az\'{e}ma--Yor martingales, stopped appropriately, are examples of
martingales which achieve equality in \eqref{eq:BD_domination}.
We can reformulate the previous result in terms of optimality of the
Az\'{e}ma--Yor stopping time, which has been studied by several authors
(Az\'{e}ma and Yor~\cite{ay2}, Gilat and Meljison~\cite{GM88}, Kertz and R\"{o}sler~\cite{KR90} and
Hobson~\cite{Hob98}).
%
\begin{coro}[(Az\'{e}ma and Yor~\cite{ay2})]
\label{cor:ay_optimality}
In the setup and notation of Theorem~\ref{th:AYor}, the distribution of
$\ovl X_{T_\psi}$ is $\mu^{\mathrm{HL}}$. In consequence, $\ovl X_{T_\psi}$
dominates stochastically the maximum of any other uniformly integrable
martingale with terminal distribution $\mu$.
\end{coro}

The result is a corollary of Theorem~\ref{thm:BD} and the fact that
the maximum~$\ovl X_{T_\psi}$ is a Hardy--Littlewood maximal r.v. associated
with\vadjust{\goodbreak}
$\mu$, which follows from Proposition~\ref{prop:q_to_U}.
Alternatively, it follows from our proof of Theorem~\ref{thm:BD} upon
observing, from the definition of $T_\psi$, that $X_{T_\psi} = w_\mu
(\ovl X_{T_\psi})$ and hence, with $P_t=X_{t\land T_\psi}$, we have
a.s. equality in \eqref{eq:bhr} for $K=w_\mu(y)$ and in consequence
$\p(\overline{X}_{T_\psi}\geq
y)=\ovl\mu^{\mathrm{HL}}(y)$.

\subsection{Optimality: a dual point of view}
We identified above $\mu^{\mathrm{HL}}$ as the maximal, relative to stochastic
order, possible distribution of maximum of a uniformly integrable
martingale with a fixed terminal law $\mu$. We   consider now the dual
problem: we look for a maximal terminal distribution of a~{uniformly}
integrable martingale with a fixed law of its maximum. We saw in~\eqref
{eq:HL_to_calls} that stochastic order of HL transforms translates into
increasing convex ordering of the underlying distributions, and we
expect the solution to the dual problem to be optimal relative to
increasing convex order.

We offer two viewpoints on this dual problem. First we will see it as
a~problem of finding the inverse operator to the Hardy--Littlewood
transform. Then we will rephrase the result in martingale terms.

Let us fix a distribution $\nu$, which the reader may think of as the
distribution of the one-sided maximum of some uniformly integrable
martingale. We look at measures $\rho$, $\int|x|\rho(dx)<\infty
$, such that $\rho^{\mathrm{HL}}$ stochastically dominates $\nu\dvtx\ovl\nu
(x)\leq\ovl\rho^{\mathrm{HL}}(x)$, $x\in\re$. We note $\mathcal{S}_\nu$
the set of such measures.
Passing to the inverses,
we can express the condition on $\rho\in\mathcal{S}_\nu$ in terms
of tail quantiles,
%
%
\begin{eqnarray}
\label{eq:def_S_set}
\rho\in\mathcal{S}_\nu \quad \Leftrightarrow \quad
\ovl q_\nu(\lambda)\leq\ovl q_{\rho^{\mathrm{HL}}}(\lambda)=\operatorname{AVaR}_{\rho}
(\lambda)=\frac{1}{\lambda}\int_0^\lambda\ovl q_\rho
(\eta)\,d\eta, \nonumber
\\[-8pt]
\\[-8pt]
\eqntext{\lambda\in[0,1],}
\end{eqnarray}
using definitions in \eqref{eq:AVaR} and \eqref{eq:HL_definition}.
Note that for existence of $\rho\in\mathcal{S}_\nu$ it is necessary that
%
%
\begin{equation}\label{eq:assume_nu}
\lambda\ovl q_\nu(\lambda)\mathop{\longrightarrow}_{\lambda\to0} 0\mbox{ which is equivalent to }
x\ovl\nu(x)\mathop{\longrightarrow}_{x\to\infty} 0 ,
\end{equation}
where the equivalence follows from the change of variables $x= \ovl
q_\nu(\lambda_x)$ and $\lambda_x \ovl q_\nu(\lambda_x)=x\ovl\nu
(x)\,dx$-a.e.

If $\nu=\mu^{\mathrm{HL}}$, then by \eqref{eq:HL_to_calls} we see instantly
that $\mu$ is the minimal element of $\mathcal{S}_\nu$ relative to
the increasing convex order. We extend this now to general measures
$\nu$.
%
\begin{theorem}\label{thm:nu_delta}
Let $\nu$ be a probability measure on $\re$. The set $\mathcal
{S}_\nu$ is nonempty if and only if $\nu$ satisfies \eqref
{eq:assume_nu}. Under \eqref{eq:assume_nu}, $\mathcal{S}_\nu$ admits
a minimal element $\nu_\Delta$ relative to the increasing convex
order, which is characterized by
\[
\lambda\ovl q_{\nu_\Delta^{\mathrm{HL}}}(\lambda)=\int_0^\lambda\ovl
q_{\nu_\Delta}(\eta)\,d\eta  \qquad \mbox{is the concave envelope of
}\lambda\ovl q_\nu(\lambda).
\]
If $\nu=\mu^{\mathrm{HL}}$ for an integrable probability measure $\mu$,
then $\nu_\Delta=\mu$.
\end{theorem}

\begin{pf}
As observed above, the case $\nu=\mu^{\mathrm{HL}}$ follows instantly from
\eqref{eq:HL_to_calls}, which in turn used the fact that $\lambda\ovl
q_{\nu}(\lambda)$ is concave and equal to $\int_0^\lambda\ovl q_\mu
(\eta)\,d\eta$.

It is natural to extend the above ideas to a general case.
Assume $\nu$ satisfies~\eqref{eq:assume_nu}, and let $G(\lambda)$ be
the concave envelope (i.e., the smallest concave majorant) of $\lambda
\ovl q_\nu(\lambda)$. If there exists a measure $\nu_\Delta$ such
that $G(\lambda)=\int_0^\lambda\ovl q_{\nu_\Delta}(\eta)\,d\eta
$, then clearly $\nu_\Delta\in\mathcal{S}_\nu$ by definition in
\eqref{eq:def_S_set}. Furthermore, since $\int_0^\lambda\ovl q_\rho
(\eta)\,d\eta$ is a concave function, we have that
%
%
\begin{equation}
\int_0^\lambda\ovl q_{\nu_\Delta}(\eta)\,d\eta\leq\int
_0^\lambda\ovl q_\rho(\eta)\,d\eta,\qquad\lambda\in[0,1],
\ \forall \rho\in\mathcal{S}_\nu .
\end{equation}
This in turn, using \eqref{eq:AVaR_to_calls}, is equivalent to $\nu
_\Delta$ being the infimum of $\rho\in\mathcal{S}_\nu$ relative to
increasing convex ordering of measures and thus being a solution to our
dual problem.

It remains to argue that $\nu_\Delta$ exists for a general $\nu$.
Recall that $-\infty\leq l_\nu<r_\nu\leq\infty$ are, respectively,
the lower and the upper bounds of the support of~$\nu$. Let $\tilde
{G}(x)$ be the (formal) Fenchel transform of $\lambda\ovl q_\nu
(\lambda)$,
%
%
\begin{equation}\label{eq:G_pseudo_Fenchel}
\tilde{G}(x)=\sup_{\lambda\in(0,1)} \bigl(\lambda\ovl q_\nu
(\lambda)-\lambda x \bigr),\qquad x\in[l_\nu,r_\nu].
\end{equation}
Observe that $\tilde{G}(x)\geq0$ thanks to assumption \eqref
{eq:assume_nu} and by definition $\tilde{G}(x)$ is convex, decreasing
and $\tilde{G}'(x)\in[-1,0]$. This implies that there exists
a~probability measure $\nu_\Delta$ such that $\tilde{G}(x)=\int
(y-x)^+\nu_\Delta(dy)=C_{\nu_\Delta}(y)$. In fact we simply have
$\ovl\nu_\Delta(x) := - \tilde{G}'(x-)$. Since $G$ was the concave
envelope of $\lambda\ovl q_\nu(\lambda)$ we can recover it as the
dual Fenchel transform of $\tilde{G}$ and, using~\eqref{eq:AVaR_2},
we have
%
%
\begin{equation}
G(\lambda)=\inf_{x\in[l_\nu,r_\nu]} \bigl(\tilde{G}(x)+x\lambda
 \bigr)=\int_0^\lambda\ovl q_{\nu_\Delta}(\eta)\,d\eta,\qquad
\lambda\in[0,1],
\end{equation}
as required. Note that we could also take $x\in\re$ above since the
infimum is always attained for $x\in[l_\nu,r_\nu]$.
\end{pf}

%
\begin{rem}
  Theorem~\ref{thm:nu_delta} synthesizes several results from Kertz
and R\"{o}sler~\cite{KertzRosler92b,KertzRosler93} as well as adds a
new interpretation of $\Delta$ operator as the inverse of $\mu\to\mu
^{\mathrm{HL}}$. Furthermore, we stress that in the proof we obtained, in fact,
a rather explicit representation which can be used to construct $\nu
_\Delta$. Namely we have $\tilde{G}(x)=C_{\nu_\Delta}(y)$ with
$\tilde{G}$ defined in \eqref{eq:G_pseudo_Fenchel}, in particular
$\ovl\nu_\Delta(x)=-\tilde{G}'(x)$. Equivalently we have $\nu
_\Delta\sim\frac{1}{\xi}G(\xi)$, for a uniform r.v. $\xi$.
\end{rem}

We offer now the second viewpoint on the problem and rephrase the
results above in martingale terms.
%
\begin{theorem}
\label{thm:HLconvexorder} Let $\nu$ be a distribution satisfying
\eqref{eq:assume_nu}, and $U(x)$ be the increasing concave envelope of
$\ovl q_\nu(1/x)$, $x\geq1$. Let $Y_t=M^U_t(N)$ for some $(N_t)$
satisfying \eqref{eq:Ndef}, $N_0=1$.\vadjust{\goodbreak}

Then $\ovl Y_\i$ dominates $\nu$ for the stochastic order and if
$(P_t)$ is any uniformly integrable martingale such that $\ovl P_\i$
dominates $\nu$ for the stochastic order, then $P_\i$ dominates $Y_\i
$ for the increasing convex order.

Furthermore, if $\nu=\mu^{\mathrm{HL}}$ and $(P_t)$ as above is
max-continuous with $P_\infty\sim\mu$ then $P$ is the Az\'
{e}ma--Yor martingale $M^{U_{\mu}}(N)$ for some $(N_t)$ satisfying
\eqref{eq:Ndef}.
\end{theorem}

\begin{pf} From \eqref{eq:assume_nu} $U$ is well defined and
observe that $U(x)=xG(1/x)$, $x\geq1$, where $G(\lambda)$ is the
concave envelope of $\lambda\ovl q_\nu(\lambda)$. From the proof of
Theorem~\ref{thm:nu_delta} we see that $\frac{1}{\lambda}G(\lambda
)=\operatorname{AVaR}_{\nu_{\Delta}}(\lambda)$ and\vspace*{1pt} hence $Y=M^{U_{\nu_\Delta
}}_t(N)$ is the Az\'{e}ma--Yor martingale associated with $\nu_\Delta
$ by Proposition~\ref{prop:q_to_U}.
Let $\mu\sim P_\infty$. By Corollary~\ref{cor:ay_optimality} the
distribution of $\ovl P_\infty$ is dominated stochastically by $\mu
^{\mathrm{HL}}$. Hence $\mu^{\mathrm{HL}}$ dominates stochastically $\nu$ and $\mu\in
\mathcal{S}_\nu$. The first part of theorem is then a corollary of
Theorem~\ref{thm:nu_delta}.

It remains to argue the last statement of the theorem. Since $P_\infty
\sim\mu$ and the distribution of $\ovl P_\infty$ dominates
stochastically $\mu^{\mathrm{HL}}$ it follows from Theorem~\ref{thm:BD} that
$\ovl P_\infty\sim\mu^{\mathrm{HL}}$. We deduce from the proof of Corollary
\ref{cor:ay_optimality} that we have an a.s. equality in \eqref
{eq:bhr} for any $y>0$ and $K=w_\mu(y)$ and hence
\[
\{P_\infty\geq w_\mu(y)\}\supseteq\{\ovl P_\infty>y\}\supseteq\{
P_\infty> w_\mu(y)\}.
\]
It follows that $P_\infty=w_\mu(\ovl P_\infty)$. Further, from
uniform integrability of $(P_t)$,
\[
\e P_\i= \e P_{\zeta_{w_\mu}^P}\leq\e w_\mu(\ovl P_{\zeta_{w_\mu
}^P})\leq\e w_\mu(\ovl P_\infty).
\]
In consequence $P_t=P_{t\land\zeta_{w_\mu}^P}$, and the statement
follows with $N_t=M^{V_\mu}_{t\land\tau^{r_\mu}(P)}(P)$ (see
Theorem~\ref{th:AYor} and Remark~\ref{rem:skoro_maxcont}).
\end{pf}

%
\begin{rem}\label{rem:U_env}
Distribution $\nu_\Delta$ can also be easily recovered from $U$
since, using Proposition~\ref{prop:q_to_U}, $h_{\nu_\Delta}(x)=\ovl
q_{\nu_\Delta}(1/x)$ is obtained from \eqref{eq:odeayequation} as
$h_{\nu_\Delta}(x)=U(x)-xU'(x)$.
\end{rem}

\subsection{Floor constraint and concave order}\label{sec:floor}
In this final section we study how Theorem~\ref{thm:HLconvexorder} can
be used to solve different
optimization problems motivated by portfolio insurance. Our insight
comes in particular from
constrained portfolio optimization problems discussed by
El Karoui and Me\-ziou~\cite{ELKM06}.

Consider $g$ an increasing function on $\re_+$ such that $\lim_{x\to
\infty}g(x)/x=0$, and let $U$ be its increasing concave envelope.
Let $N_t$ satisfy \eqref{eq:Ndef} with $N_0=1$. In the financial
context, the underlying floor is modeled by $F_t=g(N_t)$. Financial
positions can be modeled with uniformly integrable martingales, and we
are interested in choosing the optimal one, among all which dominate~$F_t$
for all $t\geq0$.
We note that it is quite remarkable that this pathwise domination
requirement turns out to be equivalent to potentially weaker conditions
of ordering of distributions.

 Finally we remark that in financial
context we often use the increasing concave order between two variables
(rather then convex). This is simply a~consequence of the fact that
utility functions are typically concave.\vspace*{-2pt}

\begin{prop} Let $F_t=g(N_t)$ be the floor process, and ${\mathcal
M}^s_F$ denote the set of uniformly integrable martingales $(P_t)$,
with $P_0=U(N_0)$ and $P_t\geq F_t$, $t\geq0$.
Then the Az\'{e}ma--Yor martingale $M^U_t(N)$ belongs to ${\mathcal
M}_F^s$ and is optimal for the concave order of the terminal values,
that is, for any increasing concave function ${G}$ and $P\in\mathcal
{M}_F^s$, $\e{G}(M^U_\infty(N))\geq\e{G}(P_\infty)$.

In fact the same result holds in the larger set ${\mathcal M}^w_F$ of
uniformly integrable martingales $(P_t)$ with $P_0=U(N_0)$ and $\p
(\ovl P_\infty\geq x)\geq\p(\ovl F_\infty\geq x)$, for all $x\in
\re$.\vspace*{-2pt}
\end{prop}

\begin{pf}
Let $\nu\sim\ovl F_\infty= g(\ovl N_\infty)$, which can also be
written as $\ovl q_\nu(\lambda)=g(\frac{1}{\lambda})$.
Note that our assumption $g(x)/x\to0$ as $x\to\i$ is equivalent to
\eqref{eq:assume_nu}. Recall that $\lambda U(\frac{1}{\lambda})$ is
the concave envelope of $\lambda g(\frac{1}{\lambda})$, $\lambda\in(0,1)$.
The result now follows from Theorem~\ref{thm:HLconvexorder}. It
suffices to note that, since $\e P_\infty= F_0=\e M^U_\infty(N)$,
increasing convex order, increasing concave order and convex order on
$P_\infty$ and $M^U_\infty(N)$ are all equivalent (cf. Shaked and
Shanthikumar~\cite{ShakedShanthikumar94}, Theorems~3.A.15 and 3.A.16).\vspace*{-2pt}
\end{pf}

If we want show the above statement only for the smaller set $\mathcal
{M}^s_F$, then we can give a direct proof as in~\cite{ELKM08}. Any
martingale $P$ which dominates~$F_t$ dominates also the smallest
supermartingale $Z_t$ which dominates $F_t$, and it is easy to see that
$Z_t=U(N_t)$.
The process $(Z_t)$ is the Snell envelope of $(g(N_t))$, as shown in
Galtchouk and Mirochnitchenko~\cite{gm} using that $U$ is an affine
function on $\{x\dvtx  U(x)>g(x)\}$.

From Proposition~\ref{prop:terminalvalues} we know that
$M_t=M^U_t(N)=\e[h(\ovl N_\infty)|\mathcal{F}_t]$, where
$h(x)=U(x)-xU'(x)$, is a uniformly integrable martingale, and we also
have $\ovl M_t=U(\ovl N_t)=\ovl Z_t$ (cf. Proposition~\ref
{prop:group}). We assume $G$ is twice
continuously differentiable, the general case following via a
limiting argument. Since $G$ is concave, $G(y)-G(x)\le G'(x)(y-x)$ for
all $x,y\ge0$.
In consequence
\begin{eqnarray*}
&&\e [G(P_\infty)-G(M_\infty) ]\\[-2pt]
&& \qquad \leq\e [G'(M_\infty
)(P_\infty-M_\infty) ]
=\e [G'(h(\ovl
 N_\infty))(P_\infty-M_\infty) ]\\[-2pt]
&& \qquad\leq\e\int_0^\infty G'(h(\ovl N_t))\,d(P_t-M_t)+
\e\int_0^\infty(P_t-M_t)G''(h(\ovl N_t))\,d(h(\ovl N_t)).
\end{eqnarray*}
The first integral is a difference of two uniformly integrable martingales
(note that $\ovl N_0>0$) and its expectation is zero. For the second
integral, recall that~$h$ is increasing and the support of $d(h(\ovl
N_t))$ is contained in the support of $d\ovl N_t$ on which $M_t=\ovl
M_t=\ovl Z_t=Z_t\leq P_t$. As $G$ is concave we see that the integral
is a.s. negative which yields the desired inequality.\vspace*{-2pt}

\section*{Acknowledgment} The authors are grateful to an anonymous
referee for careful reading of the manuscript and helpful
corrections.\vadjust{\goodbreak}

\printaddresses


\begin{thebibliography}{99}
%
%
\bibitem{Azema78}
\begin{barticle}[author]
\bauthor{\bsnm{Az{\'{e}}ma},~\bfnm{J.}\binits{J.}}
(\byear{1978}).
\btitle{Repr{\'{e}}sentation multiplicative d'une surmartingale born{\'{e}}e}.
\bjournal{Z.~Wahrsch. Verw. Gebiete}
\bvolume{45}
\bpages{191--211}.
\end{barticle}
\MR{0510025}
\endbibitem

%
%
\bibitem{AzemaYor78a}
\begin{bincollection}[author]
\bauthor{\bsnm{Az{\'{e}}ma},~\bfnm{J.}\binits{J.}} \AND
\bauthor{\bsnm{Yor},~\bfnm{M.}\binits{M.}}
(\byear{1978}).
\btitle{En guise d'introduction}.
In \bbooktitle{Temps locaux}
(\beditor{\bfnm{J.}\binits{J.}~\bsnm{Az{\'{e}}ma}} \AND
\beditor{\bfnm{M.}\binits{M.}~\bsnm{Yor}}, eds.).
\bseries{Ast\'{e}risque}
\bvolume{52\textnormal{--}53}
\bpages{3--16}.
\bpublisher{Soci\'{e}t\'{e} Math\'{e}matique de France, Paris}.
\end{bincollection}
\endbibitem

%
%
\bibitem{ay1}
\begin{bincollection}[author]
\bauthor{\bsnm{Az{\'{e}}ma},~\bfnm{Jacques}\binits{J.}} \AND
\bauthor{\bsnm{Yor},~\bfnm{Marc}\binits{M.}}
(\byear{1979}).
\btitle{Une solution simple au probl\`{e}me de {S}korokhod}.
In \bbooktitle{{S}\'{e}minaire de {P}robabilit\'{e}s, XIII}.
\bseries{Lecture Notes in Math.}
\bvolume{721}
\bpages{90--115}.
\bpublisher{Springer}, \baddress{Berlin}.
\end{bincollection}
\endbibitem

%
%
\bibitem{ay2}
\begin{bincollection}[author]
\bauthor{\bsnm{Az{\'{e}}ma},~\bfnm{Jacques}\binits{J.}} \AND
\bauthor{\bsnm{Yor},~\bfnm{Marc}\binits{M.}}
(\byear{1979}).
\btitle{{L}e probl\`{e}me de {S}korokhod: Compl\'{e}ments \`{a}
``{U}ne solution simple au probl\`{e}me de {S}korokhod.''}
In \bbooktitle{{S}\'{e}minaire de {P}robabilit\'{e}s, XIII}.
\bseries{Lecture Notes in Math.}
\bvolume{721}
\bpages{625--633}.
\bpublisher{Springer}, \baddress{Berlin}.
\end{bincollection}
\MR{0544832}
\endbibitem


%
%
\bibitem{Bach}
\begin{barticle}[author]
\bauthor{\bsnm{Bachelier},~\bfnm{Louis}\binits{L.}}
(\byear{1906}).
\btitle{Th{\'{e}}orie des probabilit{\'{e}}s continues}.
\bjournal{J. Math. Pures Appl.}
\bvolume{6}(\banumber{II})
\bpages{259--327}.
\end{barticle}
\endbibitem

%
%
\bibitem{BankElKaroui04}
\begin{barticle}[author]
\bauthor{\bsnm{Bank},~\bfnm{Peter}\binits{P.}} \AND
\bauthor{\bsnm{El~Karoui},~\bfnm{Nicole}\binits{N.}}
(\byear{2004}).
\btitle{A stochastic representation theorem with applications to optimization
and obstacle problems}.
\bjournal{Ann. Probab.}
\bvolume{32}
\bpages{1030--1067}.
\end{barticle}
\MR{2044673}
\endbibitem

%
%
\bibitem{BD63}
\begin{barticle}[author]
\bauthor{\bsnm{Blackwell},~\bfnm{David}\binits{D.}} \AND
\bauthor{\bsnm{Dubins},~\bfnm{Lester~E.}\binits{L.~E.}}
(\byear{1963}).
\btitle{A converse to the dominated convergence theorem}.
\bjournal{Illinois J. Math.}
\bvolume{7}
\bpages{508--514}.
\end{barticle}
\MR{0151572}
\endbibitem

%
%
\bibitem{BHR01}
\begin{barticle}[author]
\bauthor{\bsnm{Brown},~\bfnm{H.}\binits{H.}},
\bauthor{\bsnm{Hobson},~\bfnm{D.}\binits{D.}} \AND
\bauthor{\bsnm{Rogers},~\bfnm{L.~C.~G.}\binits{L.~C.~G.}}
(\byear{2001}).
\btitle{Robust hedging of barrier options}.
\bjournal{Math. Finance}
\bvolume{11}
\bpages{285--314}.
\end{barticle}
\MR{1839367}
\endbibitem

%
%
\bibitem{CK95}
\begin{barticle}[author]
\bauthor{\bsnm{Cvitani{\'{c}}},~\bfnm{J.}\binits{J.}} \AND
\bauthor{\bsnm{Karatzas},~\bfnm{I.}\binits{I.}}
(\byear{1995}).
\btitle{On portfolio optimization under ``drawdown'' constraints}.
\bjournal{{IMA} Lecture Notes in Mathematics {\&} Applications}
\bvolume{65}
\bpages{77--88}.
\end{barticle}
\endbibitem

%
%
\bibitem{DellacherieMeyer80}
\begin{bbook}[author]
\bauthor{\bsnm{Dellacherie},~\bfnm{Claude}\binits{C.}} \AND
\bauthor{\bsnm{Meyer},~\bfnm{Paul-Andr{\'{e}}}\binits{P.-A.}}
(\byear{1980}).
\btitle{Probabilit\'{e}s et potentiel. {C}hapitres {V} \`{a} {VIII}:
Th\'{e}orie des
martingales}, \bedition{revised} ed.
\bpublisher{Hermann}, \baddress{Paris}.
\end{bbook}
\MR{0566768}
\endbibitem

%
%
\bibitem{ELKM06}
\begin{barticle}[author]
\bauthor{\bsnm{El~Karoui},~\bfnm{Nicole}\binits{N.}} \AND
\bauthor{\bsnm{Meziou},~\bfnm{Asma}\binits{A.}}
(\byear{2006}).
\btitle{Constrained optimization with respect to stochastic dominance:
Application to portfolio insurance}.
\bjournal{Math. Finance}
\bvolume{16}
\bpages{103--117}.
\end{barticle}
\MR{2194897}
\endbibitem

%
%
\bibitem{ELKM08}
\begin{barticle}[author]
\bauthor{\bsnm{El~Karoui},~\bfnm{Nicole}\binits{N.}} \AND
\bauthor{\bsnm{Meziou},~\bfnm{Asma}\binits{A.}}
(\byear{2008}).
\btitle{Max-{P}lus decomposition of supermartingales and convex order.
{A}pplication to {A}merican options and portfolio insurance}.
\bjournal{Ann. Probab.}
\bvolume{36}
\bpages{647--697}.
\end{barticle}
\MR{2393993}
\endbibitem

%
%
\bibitem{ET08}
\begin{barticle}[author]
\bauthor{\bsnm{Elie},~\bfnm{R.}\binits{R.}} \AND
\bauthor{\bsnm{Touzi},~\bfnm{N.}\binits{N.}}
(\byear{2008}).
\btitle{Optimal lifetime consumption and investment under drawdown constraint}.
\bjournal{Finance Stoch.}
\bvolume{12}
\bpages{299--330}.
\end{barticle}
\MR{2410840}
\endbibitem

%
%
\bibitem{FS04}
\begin{bbook}[author]
\bauthor{\bsnm{F{\"{o}}llmer},~\bfnm{Hans}\binits{H.}} \AND
\bauthor{\bsnm{Schied},~\bfnm{Alexander}\binits{A.}}
(\byear{2004}).
\btitle{Stochastic Finance. An Introduction in Discrete Time},
\bedition{2nd extended} ed.
\bseries{de Gruyter Studies in Mathematics}
\bvolume{27}.
\bpublisher{de Gruyter},
\baddress{Berlin}.
\end{bbook}
\MR{2169807}
\endbibitem

%
%
\bibitem{gm}
\begin{barticle}[author]
\bauthor{\bsnm{Galtchouk},~\bfnm{L.~I.}\binits{L.~I.}} \AND
\bauthor{\bsnm{Mirochnitchenko},~\bfnm{T.~P.}\binits{T.~P.}}
(\byear{1997}).
\btitle{Optimal stopping problem for continuous local martingales and some
sharp inequalities}.
\bjournal{Stochastics Stochastics Rep.}
\bvolume{61}
\bpages{21--33}.
\end{barticle}
\MR{1473912}
\endbibitem

%
%
\bibitem{GM88}
\begin{bincollection}[author]
\bauthor{\bsnm{Gilat},~\bfnm{David}\binits{D.}} \AND
\bauthor{\bsnm{Meilijson},~\bfnm{Isaac}\binits{I.}}
(\byear{1988}).
\btitle{A simple proof of a theorem of {B}lackwell \& {D}ubins on the maximum
of a uniformly integrable martingale}.
In \bbooktitle{{S}\'{e}minaire de {P}robabilit\'{e}s, XXII}.
\bseries{Lecture Notes in Math.}
\bvolume{1321}
\bpages{214--216}.
\bpublisher{Springer}, \baddress{Berlin}.
\end{bincollection}
\MR{0960529}
\endbibitem

%
%
\bibitem{GZ93}
\begin{barticle}[author]
\bauthor{\bsnm{Grossman},~\bfnm{S.J.}\binits{S.}} \AND
\bauthor{\bsnm{Zhou},~\bfnm{Z.}\binits{Z.}}
(\byear{1993}).
\btitle{Optimal investment strategies for controlling drawdowns}.
\bjournal{Math. Finance}
\bvolume{3}
\bpages{241--276}.
\end{barticle}
\endbibitem

%
%
\bibitem{HL30}
\begin{barticle}[author]
\bauthor{\bsnm{Hardy},~\bfnm{G.~H.}\binits{G.~H.}} \AND
\bauthor{\bsnm{Littlewood},~\bfnm{J.~E.}\binits{J.~E.}}
(\byear{1930}).
\btitle{A maximal theorem with function-theoretic applications}.
\bjournal{Acta Math.}
\bvolume{54}
\bpages{81--116}.
\end{barticle}
\MR{1555303}
\endbibitem

%
%
\bibitem{Hob98}
\begin{bincollection}[author]
\bauthor{\bsnm{Hobson},~\bfnm{David~G.}\binits{D.~G.}}
(\byear{1998}).
\btitle{The maximum maximum of a martingale}.
In \bbooktitle{{S}\'{e}minaire de {P}robabilit\'{e}s, XXXII}.
\bseries{Lecture Notes in Math.}
\bvolume{1686}
\bpages{250--263}.
\bpublisher{Springer}, \baddress{Berlin}.
\end{bincollection}
\MR{1655298}
\endbibitem

%
%
\bibitem{Hu98}
\begin{barticle}[author]
\bauthor{\bsnm{Hurlimann},~\bfnm{W.}\binits{W.}}
(\byear{1998}).
\btitle{On stop-loss order and the distortion pricing principle}.
\bjournal{Astin Bull.}
\bvolume{28}
\bpages{119--134}.
\end{barticle}
\endbibitem

%
%
\bibitem{KR90}
\begin{barticle}[author]
\bauthor{\bsnm{Kertz},~\bfnm{Robert~P.}\binits{R.~P.}} \AND
\bauthor{\bsnm{R{\"{o}}sler},~\bfnm{Uwe}\binits{U.}}
(\byear{1990}).
\btitle{Martingales with given maxima and terminal distributions}.
\bjournal{Israel J. Math.}
\bvolume{69}
\bpages{173--192}.
\end{barticle}
\MR{1045372}
\endbibitem

%
%
\bibitem{KertzRosler92b}
\begin{barticle}[author]
\bauthor{\bsnm{Kertz},~\bfnm{Robert~P.}\binits{R.~P.}} \AND
\bauthor{\bsnm{R{\"{o}}sler},~\bfnm{Uwe}\binits{U.}}
(\byear{1992}).
\btitle{Stochastic and convex orders and lattices of probability
measures, with
a martingale interpretation}.
\bjournal{Israel J. Math.}
\bvolume{77}
\bpages{129--164}.
\end{barticle}
\MR{1194790}
\endbibitem

%
%
\bibitem{KertzRosler93}
\begin{bincollection}[author]
\bauthor{\bsnm{Kertz},~\bfnm{Robert~P.}\binits{R.~P.}} \AND
\bauthor{\bsnm{R{\"{o}}sler},~\bfnm{Uwe}\binits{U.}}
(\byear{1993}).
\btitle{Hyperbolic-concave functions and {H}ardy--{L}ittlewood maximal
functions}.
In \bbooktitle{Stochastic Inequalities ({S}eattle, {WA}, 1991)}.
\bseries{IMS Lecture Notes Monogr. Ser.}
\bvolume{22}
\bpages{196--210}.
\bpublisher{IMS}, \baddress{Hayward, CA}.
\end{bincollection}
\endbibitem

%
%
\bibitem{MRY08}
\begin{barticle}[author]
\bauthor{\bsnm{Madan},~\bfnm{Dilip}\binits{D.}},
\bauthor{\bsnm{Roynette},~\bfnm{Bernard}\binits{B.}} \AND
\bauthor{\bsnm{Yor},~\bfnm{Marc}\binits{M.}}
(\byear{2008}).
\btitle{Option prices as probabilities}.
\bjournal{Finance Research Letters}
\bvolume{5}
\bpages{79--87}.
\end{barticle}
\endbibitem

%
%
\bibitem{NY06}
\begin{barticle}[author]
\bauthor{\bsnm{Nikeghbali},~\bfnm{Ashkan}\binits{A.}} \AND
\bauthor{\bsnm{Yor},~\bfnm{Marc}\binits{M.}}
(\byear{2006}).
\btitle{Doob's maximal identity, multiplicative decompositions and enlargements
of filtrations}.
\bjournal{Illinois J. Math.}
\bvolume{50}
\bpages{791--814} (electronic).
\end{barticle}
\MR{2247846}
\endbibitem

%
%
\bibitem{genealogia}
\begin{barticle}[author]
\bauthor{\bsnm{Ob{\l}{\'{o}}j},~\bfnm{Jan}\binits{J.}}
(\byear{2004}).
\btitle{The {S}korokhod embedding problem and its offspring}.
\bjournal{Probab. Surv.}
\bvolume{1}
\bpages{321--392}.
\end{barticle}
\MR{2068476}
\endbibitem

%
%
\bibitem{jamartcarac}
\begin{barticle}[author]
\bauthor{\bsnm{Ob{\l}{\'{o}}j},~\bfnm{Jan}\binits{J.}}
(\byear{2006}).
\btitle{A complete characterization of local martingales which are
functions of
{B}rownian motion and its supremum}.
\bjournal{Bernoulli}
\bvolume{12}
\bpages{955--969}.
\end{barticle}
\MR{2274851}
\endbibitem

%
%
\bibitem{oy}
\begin{bincollection}[author]
\bauthor{\bsnm{Ob{\l}{\'{o}}j},~\bfnm{Jan}\binits{J.}} \AND
\bauthor{\bsnm{Yor},~\bfnm{Marc}\binits{M.}}
(\byear{2006}).
\btitle{On local martingale and its supremum: Harmonic functions and beyond}.
In \bbooktitle{From Stochastic Calculus to Mathematical Finance}
(\beditor{\bfnm{Yuri}\binits{Y.}~\bsnm{Kabanov}},
\beditor{\bfnm{Robert}\binits{R.}~\bsnm{Lipster}} \AND
\beditor{\bfnm{Jordan}\binits{J.}~\bsnm{Stoyanov}}, eds.)
\bpages{517--534}.
\bpublisher{Springer},
\baddress{Berlin}.
\end{bincollection}
\MR{2234288}
\endbibitem

%
%
\bibitem{PRY10}
\begin{bbook}[author]
\bauthor{\bsnm{Profeta},~\bfnm{C.}\binits{C.}},
\bauthor{\bsnm{Roynette},~\bfnm{B.}\binits{B.}} \AND
\bauthor{\bsnm{Yor},~\bfnm{M.}\binits{M.}}
(\byear{2010}).
\btitle{Option Prices as Probabilities}.
\bpublisher{Springer},
\baddress{Berlin}.
\end{bbook}
\MR{2582990}
\endbibitem

%
%
\bibitem{ry}
\begin{bbook}[author]
\bauthor{\bsnm{Revuz},~\bfnm{Daniel}\binits{D.}} \AND
\bauthor{\bsnm{Yor},~\bfnm{Marc}\binits{M.}}
(\byear{1999}).
\btitle{Continuous Martingales and {B}rownian Motion}.
\bseries{Grundlehren der Mathematischen Wissenschaften [Fundamental Principles
of Mathematical Sciences]}
\bvolume{293}.
\bpublisher{Springer}, \baddress{Berlin}.
\end{bbook}
\MR{1725357}
\endbibitem

%
%
\bibitem{RVY06}
\begin{barticle}[author]
\bauthor{\bsnm{Roynette},~\bfnm{Bernard}\binits{B.}},
\bauthor{\bsnm{Vallois},~\bfnm{Pierre}\binits{P.}} \AND
\bauthor{\bsnm{Yor},~\bfnm{Marc}\binits{M.}}
(\byear{2006}).
\btitle{Limiting laws associated with {B}rownian motion perturbated by its
maximum, minimum and local time {II}}.
\bjournal{Studia Sci. Math. Hungar.}
\bvolume{43}
\bpages{295--360}.
\end{barticle}
\MR{2253307}
\endbibitem

%
%
\bibitem{ShakedShanthikumar94}
\begin{bbook}[author]
\bauthor{\bsnm{Shaked},~\bfnm{Moshe}\binits{M.}} \AND
\bauthor{\bsnm{Shanthikumar},~\bfnm{J.~George}\binits{J.~G.}}
(\byear{1994}).
\btitle{Stochastic Orders and Their Applications}.
\bpublisher{Academic Press}, \baddress{Boston, MA}.
\end{bbook}
\MR{1278322}
\endbibitem
\end{thebibliography}
\end{document}